
\magnification1200
\baselineskip15pt
\normalbaselineskip=15pt

\newread\AUX\immediate\openin\AUX=\jobname.aux
\newcount\relFnno
\def\ref#1{\expandafter\edef\csname#1\endcsname}
\ifeof\AUX\immediate\write16{\jobname.aux gibt es nicht!}\else
\input \jobname.aux
\fi\immediate\closein\AUX



\def\ignore{\bgroup
\catcode`\;=0\catcode`\^^I=14\catcode`\^^J=14\catcode`\^^M=14
\catcode`\ =14\catcode`\!=14\catcode`\"=14\catcode`\#=14\catcode`\$=14
\catcode`\&=14\catcode`\'=14\catcode`\(=14\catcode`\)=14\catcode`\*=14
\catcode`+=14\catcode`\,=14\catcode`\-=14\catcode`\.=14\catcode`\/=14
\catcode`\0=14\catcode`\1=14\catcode`\2=14\catcode`\3=14\catcode`\4=14
\catcode`\5=14\catcode`\6=14\catcode`\7=14\catcode`\8=14\catcode`\9=14
\catcode`\:=14\catcode`\<=14\catcode`\==14\catcode`\>=14\catcode`\?=14
\catcode`\@=14\catcode`\A=14\catcode`\B=14\catcode`\C=14\catcode`\D=14
\catcode`\E=14\catcode`\F=14\catcode`\G=14\catcode`\H=14\catcode`\I=14
\catcode`\J=14\catcode`\K=14\catcode`\L=14\catcode`\M=14\catcode`\N=14
\catcode`\O=14\catcode`\P=14\catcode`\Q=14\catcode`\R=14\catcode`\S=14
\catcode`\T=14\catcode`\U=14\catcode`\V=14\catcode`\W=14\catcode`\X=14
\catcode`\Y=14\catcode`\Z=14\catcode`\[=14\catcode`\\=14\catcode`\]=14
\catcode`\^=14\catcode`\_=14\catcode`\`=14\catcode`\a=14\catcode`\b=14
\catcode`\c=14\catcode`\d=14\catcode`\e=14\catcode`\f=14\catcode`\g=14
\catcode`\h=14\catcode`\i=14\catcode`\j=14\catcode`\k=14\catcode`\l=14
\catcode`\m=14\catcode`\n=14\catcode`\o=14\catcode`\p=14\catcode`\q=14
\catcode`\r=14\catcode`\s=14\catcode`\t=14\catcode`\u=14\catcode`\v=14
\catcode`\w=14\catcode`\x=14\catcode`\y=14\catcode`\z=14\catcode`\{=14
\catcode`\|=14\catcode`\}=14\catcode`\~=14\catcode`\^^?=14
\Ignoriere}
\def\Ignoriere#1\;{\egroup}

\newcount\itemcount
\def\resetitem{\global\itemcount0}\resetitem
\newcount\Itemcount
\Itemcount0
\newcount\maxItemcount
\maxItemcount=0

\def\FILTER\fam\itfam\tenit#1){#1}

\def\Item#1{\global\advance\itemcount1
\edef\TEXT{{\it\romannumeral\itemcount)}}%
\ifx?#1?\relax\else
\ifnum#1>\maxItemcount\global\maxItemcount=#1\fi
\expandafter\ifx\csname I#1\endcsname\relax\else
\edef\testA{\csname I#1\endcsname}
\expandafter\expandafter\def\expandafter\testB\testA
\edef\testC{\expandafter\FILTER\testB}
\edef\testD{\csname0\testC0\endcsname}\fi
\edef\testE{\csname0\romannumeral\itemcount0\endcsname}
\ifx\testD\testE\relax\else
\immediate\write16{I#1 hat sich geaendert!}\fi
\expandwrite\AUX{\neverexpand\ref{I#1}{\TEXT}}\fi
\item{\TEXT}}

\def\today{\number\day.~\ifcase\month\or
  Januar\or Februar\or M{\"a}rz\or April\or Mai\or Juni\or
  Juli\or August\or September\or Oktober\or November\or Dezember\fi
  \space\number\year}
\font\sevenex=cmex7
\scriptfont3=\sevenex
\font\fiveex=cmex10 scaled 500
\scriptscriptfont3=\fiveex
\def\A{{\bf A}}
\def\G{{\bf G}}
\def\P{{\bf P}}

\def\phi{\varphi}
\def\epsilon{\varepsilon}

\def\theta{\vartheta}
\def\uauf{\lower1.7pt\hbox to 3pt{%
\vbox{\offinterlineskip
\hbox{\vbox to 8.5pt{\leaders\vrule width0.2pt\vfill}%
\kern-.3pt\hbox{\lams\char"76}\kern-0.3pt%
$\raise1pt\hbox{\lams\char"76}$}}\hfil}}

\def\title#1{\par
{\baselineskip1.5\baselineskip\rightskip0pt plus 5truecm
\leavevmode\vskip0truecm\noindent\font\BF=cmbx10 scaled \magstep2\BF #1\par}
\vskip1truecm
\leftline{\font\CSC=cmcsc10
{\CSC Friedrich Knop}}
\leftline{Department of Mathematics, Rutgers University, New Brunswick NJ
08903, USA}
\leftline{knop@math.rutgers.edu}
\vskip1truecm
\par}

\def\cite#1{\expandafter\ifx\csname#1\endcsname\relax
{\bf?}\immediate\write16{#1 ist nicht definiert!}\else\csname#1\endcsname\fi}
\def\expandwrite#1#2{\edef\next{\write#1{#2}}\next}
\def\neverexpand{\noexpand\noexpand\noexpand}
\def\strip#1\ {}
\def\ncite#1{\expandafter\ifx\csname#1\endcsname\relax
{\bf?}\immediate\write16{#1 ist nicht definiert!}\else
\expandafter\expandafter\expandafter\strip\csname#1\endcsname\fi}
\newwrite\AUX
\immediate\openout\AUX=\jobname.aux
\font\eightrm=cmr8\font\sixrm=cmr6
\font\eighti=cmmi8
\font\eightit=cmti8
\font\eightbf=cmbx8
\font\eightcsc=cmcsc10 scaled 833
\def\eightpoint{%
\textfont0=\eightrm\scriptfont0=\sixrm\def\rm{\fam0\eightrm}%
\textfont1=\eighti
\textfont\bffam=\eightbf\def\bf{\fam\bffam\eightbf}%
\textfont\itfam=\eightit\def\it{\fam\itfam\eightit}%
\def\csc{\eightcsc}%
\setbox\strutbox=\hbox{\vrule height7pt depth2pt width0pt}%
\normalbaselineskip=0,8\normalbaselineskip\normalbaselines\rm}
\newcount\absFnno\absFnno1
\write\AUX{\relFnno1}
\newif\ifMARKE\MARKEtrue
{\catcode`\@=11
\gdef\footnote{\ifMARKE\edef\@sf{\spacefactor\the\spacefactor}\/%
$^{\cite{Fn\the\absFnno}}$\@sf\fi
\MARKEtrue
\insert\footins\bgroup\eightpoint
\interlinepenalty100\let\par=\endgraf
\leftskip=0pt\rightskip=0pt
\splittopskip=10pt plus 1pt minus 1pt \floatingpenalty=20000\smallskip
\item{$^{\cite{Fn\the\absFnno}}$}%
\expandwrite\AUX{\neverexpand\ref{Fn\the\absFnno}{\neverexpand\the\relFnno}}%
\global\advance\absFnno1\write\AUX{\advance\relFnno1}%
\bgroup\strut\aftergroup\@foot\let\next}}
\skip\footins=12pt plus 2pt minus 4pt
\dimen\footins=30pc
\output={\plainoutput\immediate\write\AUX{\relFnno1}}
\newcount\Abschnitt\Abschnitt0
\def\beginsection#1. #2 \par{\advance\Abschnitt1%
\vskip0pt plus.10\vsize\penalty-250
\vskip0pt plus-.10\vsize\bigskip\vskip\parskip
\edef\TEST{\number\Abschnitt}
\expandafter\ifx\csname#1\endcsname\TEST\relax\else
\immediate\write16{#1 hat sich geaendert!}\fi
\expandwrite\AUX{\neverexpand\ref{#1}{\TEST}}
\leftline{\marginnote{#1}\bf\number\Abschnitt. \ignorespaces#2}%
\nobreak\smallskip\noindent\SATZ1\GNo0}
\def\Proof:{\par\noindent{\it Proof:}}
\def\Remark:{\ifdim\lastskip<\medskipamount\removelastskip\medskip\fi
\noindent{\bf Remark:}}
\def\Remarks:{\ifdim\lastskip<\medskipamount\removelastskip\medskip\fi
\noindent{\bf Remarks:}}
\def\Definition:{\ifdim\lastskip<\medskipamount\removelastskip\medskip\fi
\noindent{\bf Definition:}}
\def\Example:{\ifdim\lastskip<\medskipamount\removelastskip\medskip\fi
\noindent{\bf Example:}}
\def\Examples:{\ifdim\lastskip<\medskipamount\removelastskip\medskip\fi
\noindent{\bf Examples:}}
\newif\ifmarginalnotes\marginalnotesfalse
\newif\ifmarginalwarnings\marginalwarningstrue

\def\marginnote#1{\ifmarginalnotes\hbox to 0pt{\eightpoint\hss #1\ }\fi}

\def\strutdepth{\dp\strutbox}
\def\Randbem#1#2{\ifmarginalwarnings
{#1}\strut
\setbox0=\vtop{\eightpoint
\rightskip=0pt plus 6mm\hfuzz=3pt\hsize=16mm\noindent\leavevmode#2}%
\vadjust{\kern-\strutdepth
\vtop to \strutdepth{\kern-\ht0
\hbox to \hsize{\kern-16mm\kern-6pt\box0\kern6pt\hfill}\vss}}\fi}

\def\Zitat!{\Randbem{\bf?}{\bf Zitat}}

\newcount\SATZ\SATZ1
\def\proclaim #1. #2\par{\ifdim\lastskip<\medskipamount\removelastskip
\medskip\fi
\noindent{\bf#1.\ }{\it#2}\Par
\ifdim\lastskip<\medskipamount\removelastskip\goodbreak\medskip\fi}
\def\Aussage#1{\expandafter\def\csname#1\endcsname##1.{\resetitem
\ifx?##1?\relax\else
\edef\TEST{#1\penalty10000\ \number\Abschnitt.\number\SATZ}
\expandafter\ifx\csname##1\endcsname\TEST\relax\else
\immediate\write16{##1 hat sich geaendert!}\fi
\expandwrite\AUX{\neverexpand\ref{##1}{\TEST}}\fi
\proclaim {\marginnote{##1}\number\Abschnitt.\number\SATZ. #1\global\advance\SATZ1}.}}
\Aussage{Theorem}
\Aussage{Proposition}
\Aussage{Corollary}
\Aussage{Lemma}
\font\la=lasy10
\def\strich{\hbox{$\vcenter{\hbox
to 1pt{\leaders\hrule height -0,2pt depth 0,6pt\hfil}}$}}
\def\dashedrightarrow{\hbox{%
\hbox to 0,5cm{\leaders\hbox to 2pt{\hfil\strich\hfil}\hfil}%
\kern-2pt\hbox{\la\char\string"29}}}

\def\Bindestrich{\penalty10000-\hskip0pt}
\let\_=\Bindestrich
\def\.{{\sfcode`.=1000.}}

\def\Rechts#1{\rlap{$\scriptstyle#1$}}
\def\Par{\par}
\def\:={\mathrel{\raise0,9pt\hbox{.}\kern-2,77779pt
\raise3pt\hbox{.}\kern-2,5pt=}}
\def\=:{\mathrel{=\kern-2,5pt\raise0,9pt\hbox{.}\kern-2,77779pt
\raise3pt\hbox{.}}} 
\def\into{\hookrightarrow}
\def\pfeil{\rightarrow}
\def\untenPf{\downarrow}

\def\pf#1{\buildrel#1\over\rightarrow}
\def\Pf#1{\buildrel#1\over\longrightarrow}

\def\Ugleich{\hbox{$\cup$\kern.5pt\vrule depth -0.5pt}}
\def\|#1|{\mathop{\rm#1}\nolimits}
\def\<{\langle}
\def\>{\rangle}
\let\Times=\times
\def\times{\mathop{\Times}}
\let\Otimes=\otimes
\def\otimes{\mathop{\Otimes}}
\catcode`\@=11
\def\hex#1{\ifcase#1 0\or1\or2\or3\or4\or5\or6\or7\or8\or9\or A\or B\or
C\or D\or E\or F\else\message{Warnung: Setze hex#1=0}0\fi}
\def\fontdef#1:#2,#3,#4.{%
\alloc@8\fam\chardef\sixt@@n\FAM
\ifx!#2!\else\expandafter\font\csname text#1\endcsname=#2
\textfont\the\FAM=\csname text#1\endcsname\fi
\ifx!#3!\else\expandafter\font\csname script#1\endcsname=#3
\scriptfont\the\FAM=\csname script#1\endcsname\fi
\ifx!#4!\else\expandafter\font\csname scriptscript#1\endcsname=#4
\scriptscriptfont\the\FAM=\csname scriptscript#1\endcsname\fi
\expandafter\edef\csname #1\endcsname{\fam\the\FAM\csname text#1\endcsname}
\expandafter\edef\csname hex#1fam\endcsname{\hex\FAM}}
\catcode`\@=12 

\fontdef Ss:cmss10,,.
\fontdef Fr:eufm10,eufm7,eufm5.


%
\fontdef bbb:msbm10,msbm7,msbm5.
\fontdef mbf:cmmib10,cmmib7,.
\fontdef msa:msam10,msam7,msam5.

\def\NN{{\bbb N}}
\def\QQ{{\bbb Q}}

\def\ZZ{{\bbb Z}}

\def\cE{{\cal E}}\def\cH{{\cal H}}
\def\cJ{{\cal J}}\def\cL{{\cal L}}
\def\cM{{\cal M}}\def\cP{{\cal P}}

\mathchardef\leer=\string"0\hexbbbfam3F
\mathchardef\subsetneq=\string"3\hexbbbfam24
\mathchardef\semidir=\string"2\hexbbbfam6E
\mathchardef\dirsemi=\string"2\hexbbbfam6F
\mathchardef\haken=\string"2\hexmsafam78
\mathchardef\auf=\string"3\hexmsafam10
\let\OL=\overline
\def\overline#1{{\hskip1pt\OL{\hskip-1pt#1\hskip-.3pt}\hskip.3pt}}
\def\Aq{{\overline{A}}}


\def\Hq{{\overline{H}}}

\def\oq{{\overline{o}}}
\def\pq{{\overline{p}}}

\def\vq{{\overline{v}}}

%
\newdimen\Parindent
\Parindent=\parindent


\abovedisplayskip 9.0pt plus 3.0pt minus 3.0pt
\belowdisplayskip 9.0pt plus 3.0pt minus 3.0pt
\newdimen\Grenze\Grenze2\Parindent\advance\Grenze1em
\newdimen\Breite
\newbox\DpBox
\def\NewDisplay#1
#2$${\Breite\hsize\advance\Breite-\hangindent
\setbox\DpBox=\hbox{\hskip2\Parindent$\displaystyle{%
\ifx0#1\relax\else\eqno{#1}\fi#2}$}%
\ifnum\predisplaysize<\Grenze\abovedisplayskip\abovedisplayshortskip
\belowdisplayskip\belowdisplayshortskip\fi
\global\futurelet\nexttok\WEITER}
\def\WEITER{\ifx\nexttok\qed\expandafter\leftQEDdisplay
\else\leftdisplay\fi}
\def\leftdisplay{\hskip-\hangindent\leftline{\box\DpBox}$$}
\def\leftQEDdisplay{\hskip-\hangindent
\line{\copy\DpBox\hfill\lower\dp\DpBox\copy\QEDbox}%
\belowdisplayskip0pt$$\bigskip\let\nexttok=}
\everydisplay{\NewDisplay}
\newcount\GNo\GNo=0
\newcount\maxEqNo\maxEqNo=0
\def\eqno#1{%
\global\advance\GNo1
\edef\FTEST{$(\number\Abschnitt.\number\GNo)$}
\ifx?#1?\relax\else
\ifnum#1>\maxEqNo\global\maxEqNo=#1\fi%
\expandafter\ifx\csname E#1\endcsname\FTEST\relax\else
\immediate\write16{E#1 hat sich geaendert!}\fi
\expandwrite\AUX{\neverexpand\ref{E#1}{\FTEST}}\fi
\llap{\hbox to 40pt{\marginnote{#1}\FTEST\hfill}}}

\catcode`@=11
\def\eqalignno#1{\null\!\!\vcenter{\openup\jot\m@th\ialign{\eqno{##}\hfil
&\strut\hfil$\displaystyle{##}$&$\displaystyle{{}##}$\hfil\crcr#1\crcr}}\,}
\catcode`@=12

\newbox\QEDbox
\newbox\nichts\setbox\nichts=\vbox{}\wd\nichts=2mm\ht\nichts=2mm
\setbox\QEDbox=\hbox{\vrule\vbox{\hrule\copy\nichts\hrule}\vrule}
\def\qed{\leavevmode\unskip\hfil\null\nobreak\hfill\copy\QEDbox\medbreak}
\newdimen\HIindent
\newbox\HIbox
\def\setHI#1{\setbox\HIbox=\hbox{#1}\HIindent=\wd\HIbox}
\def\HI#1{\par\hangindent\HIindent\hangafter=0\noindent\leavevmode
\llap{\hbox to\HIindent{#1\hfil}}\ignorespaces}

\newdimen\maxSpalbr
\newdimen\altSpalbr
\newcount\Zaehler


\newif\ifxxx

{\catcode`/=\active

\gdef\beginrefs{%
\xxxfalse
\catcode`/=\active
\def/{\string/\ifxxx\hskip0pt\fi}
\def\TText##1{{\xxxtrue\tt##1}}
\expandafter\ifx\csname Spaltenbreite\endcsname\relax
\def\Spaltenbreite{1cm}\immediate\write16{Spaltenbreite undefiniert!}\fi
\expandafter\altSpalbr\Spaltenbreite
\maxSpalbr0pt
\gdef\alt{}
\def\\##1\relax{%
\gdef\neu{##1}\ifx\alt\neu\global\advance\Zaehler1\else
\xdef\alt{\neu}\global\Zaehler=1\fi\xdef\SigText{##1\the\Zaehler}}
\def\L|Abk:##1|Sig:##2|Au:##3|Tit:##4|Zs:##5|Bd:##6|S:##7|J:##8|xxx:##9||{%
\def\SigText{##2}\global\setbox0=\hbox{##2\relax}
\edef\TEST{[\SigText]}
\expandafter\ifx\csname##1\endcsname\TEST\relax\else
\immediate\write16{##1 hat sich geaendert!}\fi
\expandwrite\AUX{\neverexpand\ref{##1}{\TEST}}
\setHI{[\SigText]\ }
\ifnum\HIindent>\maxSpalbr\maxSpalbr\HIindent\fi
\ifnum\HIindent<\altSpalbr\HIindent\altSpalbr\fi
\HI{\marginnote{##1}[\SigText]}
\ifx-##3\relax\else{##3}: \fi
\ifx-##4\relax\else{##4}{\sfcode`.=3000.} \fi
\ifx-##5\relax\else{\it ##5\/} \fi
\ifx-##6\relax\else{\bf ##6} \fi
\ifx-##8\relax\else({##8})\fi
\ifx-##7\relax\else, {##7}\fi
\ifx-##9\relax\else, \TText{##9}\fi\Par}
\def\B|Abk:##1|Sig:##2|Au:##3|Tit:##4|Reihe:##5|Verlag:##6|Ort:##7|J:##8|xxx:##9||{%
\def\SigText{##2}\global\setbox0=\hbox{##2\relax}
\edef\TEST{[\SigText]}
\expandafter\ifx\csname##1\endcsname\TEST\relax\else
\immediate\write16{##1 hat sich geaendert!}\fi
\expandwrite\AUX{\neverexpand\ref{##1}{\TEST}}
\setHI{[\SigText]\ }
\ifnum\HIindent>\maxSpalbr\maxSpalbr\HIindent\fi
\ifnum\HIindent<\altSpalbr\HIindent\altSpalbr\fi
\HI{\marginnote{##1}[\SigText]}
\ifx-##3\relax\else{##3}: \fi
\ifx-##4\relax\else{##4}{\sfcode`.=3000.} \fi
\ifx-##5\relax\else{(##5)} \fi
\ifx-##7\relax\else{##7:} \fi
\ifx-##6\relax\else{##6}\fi
\ifx-##8\relax\else{ ##8}\fi
\ifx-##9\relax\else, \TText{##9}\fi\Par}
\def\Pr|Abk:##1|Sig:##2|Au:##3|Artikel:##4|Titel:##5|Hgr:##6|Reihe:{%
\def\SigText{##2}\global\setbox0=\hbox{##2\relax}
\edef\TEST{[\SigText]}
\expandafter\ifx\csname##1\endcsname\TEST\relax\else
\immediate\write16{##1 hat sich geaendert!}\fi
\expandwrite\AUX{\neverexpand\ref{##1}{\TEST}}
\setHI{[\SigText]\ }
\ifnum\HIindent>\maxSpalbr\maxSpalbr\HIindent\fi
\ifnum\HIindent<\altSpalbr\HIindent\altSpalbr\fi
\HI{\marginnote{##1}[\SigText]}
\ifx-##3\relax\else{##3}: \fi
\ifx-##4\relax\else{##4}{\sfcode`.=3000.} \fi
\ifx-##5\relax\else{In: \it ##5}. \fi
\ifx-##6\relax\else{(##6)} \fi\PrII}
\def\PrII##1|Bd:##2|Verlag:##3|Ort:##4|S:##5|J:##6|xxx:##7||{%
\ifx-##1\relax\else{##1} \fi
\ifx-##2\relax\else{\bf ##2}, \fi
\ifx-##4\relax\else{##4:} \fi
\ifx-##3\relax\else{##3} \fi
\ifx-##6\relax\else{##6}\fi
\ifx-##5\relax\else{, ##5}\fi
\ifx-##7\relax\else, \TText{##7}\fi\Par}
\bgroup
\baselineskip12pt
\parskip2.5pt plus 1pt
\hyphenation{Hei-del-berg Sprin-ger}
\sfcode`.=1000
\beginsection References. References

}}

\def\endrefs{%
\expandwrite\AUX{\neverexpand\ref{Spaltenbreite}{\the\maxSpalbr}}
\ifnum\maxSpalbr=\altSpalbr\relax\else
\immediate\write16{Spaltenbreite hat sich geaendert!}\fi
\egroup\write16{Letzte Gleichung: E\the\maxEqNo}
\write16{Letzte Aufzaehlung: I\the\maxItemcount}}


\catcode`\@11\relax
\newif\ify@autoscale \y@autoscaletrue \def\Yautoscale#1{\ifnum #1=0
  \y@autoscalefalse\else\y@autoscaletrue\fi}
\newdimen\y@b@xdim
\newdimen\y@boxdim \y@boxdim=13pt
\def\Yboxdim#1{\y@autoscalefalse\y@boxdim=#1}
\newdimen\y@linethick    \y@linethick=.3pt
\def\Ylinethick#1{\y@linethick=#1}
\newskip\y@interspace \y@interspace=0ex plus 0.3ex
\def\Yinterspace#1{\y@interspace=#1}
\newif\ify@vcenter   \y@vcenterfalse
\def\Yvcentermath#1{\ifnum #1=0 \y@vcenterfalse\else\y@vcentertrue\fi}
\newif\ify@stdtext   \y@stdtextfalse
\def\Ystdtext#1{\ifnum #1=0 \y@stdtextfalse\else\y@stdtexttrue\fi}
\newif\ify@enable@skew   \y@enable@skewfalse
\expandafter\ifx\csname enableskew\endcsname\relax
 \y@enable@skewfalse \else \y@enable@skewtrue\fi
\def\y@vr{\vrule height0.8\y@b@xdim width\y@linethick depth 0.2\y@b@xdim}
\def\y@emptybox{\y@vr\hbox to \y@b@xdim{\hfil}}
\ify@enable@skew
 \def\y@abcbox#1{\if :#1\else
   \y@vr\hbox to \y@b@xdim{\hfil#1\hfil}\fi}
 \def\y@mathabcbox#1{\if :#1\else
   \y@vr\hbox to \y@b@xdim{\hfil$#1$\hfil}\fi}
\else
 \def\y@abcbox#1{\y@vr\hbox to \y@b@xdim{\hfil#1\hfil}}
 \def\y@mathabcbox#1{\y@vr\hbox to \y@b@xdim{\hfil$#1$\hfil}}
\fi
\def\y@setdim{%
  \ify@autoscale%
   \ifvoid1\else\typeout{Package youngtab: box1 not free! Expect an
     error!}\fi%
   \setbox1=\hbox{A}\y@b@xdim=1.6\ht1 \setbox1=\hbox{}\box1%
  \else\y@b@xdim=\y@boxdim \advance\y@b@xdim by -2\y@linethick
  \fi}
\newcount\y@counter
\newif\ify@islastarg
\def\y@lastargtest#1,#2 {\if\space #2 \y@islastargtrue
  \else\y@islastargfalse\fi}
\def\y@emptyboxes#1{\y@counter=#1\loop\ifnum\y@counter>0
  \advance\y@counter by -1 \y@emptybox\repeat}
\def\y@nelineemptyboxes#1{%
  \vbox{%
    \hrule height\y@linethick%
    \hbox{\y@emptyboxes{#1}\y@vr}
    \hrule height\y@linethick}\vskip-\y@linethick}
\def\yng(#1){%
  \y@setdim%
  \hskip\y@interspace%
  \ifmmode\ify@vcenter\vcenter\fi\fi{%
  \y@lastargtest#1,
  \vbox{\offinterlineskip
    \ify@islastarg
     \y@nelineemptyboxes{#1}
    \else
     \y@ungempty(#1)
    \fi}}\hskip\y@interspace}
\def\y@ungempty(#1,#2){%
  \y@nelineemptyboxes{#1}
  \y@lastargtest#2,
  \ify@islastarg
   \y@nelineemptyboxes{#2}
  \else
   \y@ungempty(#2)
  \fi}
\def\y@nelettertest#1#2. {\if\space #2 \y@islastargtrue
  \else\y@islastargfalse\fi}
\def\y@abcboxes#1#2.{%
  \ify@stdtext\y@abcbox#1\else\y@mathabcbox#1\fi%
  \y@nelettertest #2.
  \ify@islastarg\unskip%
   \ify@stdtext\y@abcbox{#2}\else\y@mathabcbox{#2}\fi%
  \else\y@abcboxes#2.\fi}
 \newdimen\y@full@b@xdim
 \newcount\y@m@veright@cnt
\ify@enable@skew
 \def\y@get@m@veright@cnt#1#2.{%
   \if :#1 \advance\y@m@veright@cnt by 1\y@get@m@veright@cnt#2.\fi}
 \let\y@setdim@=\y@setdim
 \def\y@setdim{%
   \y@setdim@ \y@full@b@xdim=\y@b@xdim
   \advance\y@full@b@xdim by 1\y@linethick}
 \def\y@m@veright@ifskew#1{
   \y@m@veright@cnt=0 \y@get@m@veright@cnt#1.
   \moveright \y@m@veright@cnt\y@full@b@xdim}
\else
 \def\y@m@veright@ifskew#1{}
\fi
\def\y@nelineabcboxes#1{%
\if\space #1 \vbox{\hbox{\y@vr}}\else
%
  \y@nelettertest #1.
  \ify@islastarg
   \y@m@veright@ifskew{#1}
    \vbox{
      \hrule height\y@linethick%
      \hbox{\ify@stdtext\y@abcbox#1\else\y@mathabcbox#1\fi\y@vr}
      \hrule height\y@linethick}\vskip-\y@linethick
  \else
   \y@m@veright@ifskew{#1}
    \vbox{
      \hrule height\y@linethick%
      \hbox{\y@abcboxes #1.\y@vr}%
      \hrule height\y@linethick}\vskip-\y@linethick
  \fi\fi}
\def\young(#1){%
  \y@setdim%
  \hskip\y@interspace%
  \y@lastargtest#1,
  \ifmmode\ify@vcenter\vcenter\fi\fi{%
  \vbox{\offinterlineskip
    \ify@islastarg\y@nelineabcboxes{#1}%
    \else\y@ungabc(#1)%
    \fi}}\hskip\y@interspace}
\def\y@ungabc(#1,#2){%
  \y@nelineabcboxes{#1}%
  \y@lastargtest#2,
  \ify@islastarg\y@nelineabcboxes{#2}%
  \else\y@ungabc(#2)%
  \fi}
\catcode`\@12\relax

\Aussage{Conjecture}

\def\v{^\vee}

\def\Wf{W_{\mskip-5mu f}}

\def\uH{{\underline H}}
\def\uM{{\underline M}}

\def\oq{\omega}
\def\Mp{\cM^{\rm pol}}
\def\Hp{{\cH^{\rm pol}}}
\def\Wp{{W^{\rm pol}}}
\def\Ma{\cM^{\rm as}}
\def\Pa{\cP^{\rm as}}
\def\npd{{\rm pl}}
\def\plim{\mathop{\vtop{\hbox{\rm lim}\vskip-8pt\hbox{$\longleftarrow$}}}}
\fontdef bsy:cmmib10,cmmib7,cmmib5.

\mathchardef\blambda=\string"0\hexbsyfam15
\mathchardef\bmu=\string"0\hexbsyfam16

\def\punkt{\hbox{\vrule height9pt width 10.95pt depth 2.2pt}}
\def\leer{\hbox{}}


\title{Composition Kostka functions}

{\narrower\noindent Macdonald defined two\_parameter Kostka functions
$K_{\lambda\mu}(q,t)$ where $\lambda$, $\mu$ are partitions. The main
purpose of this paper is to extend his definition to include all
compositions as indices. Following Macdonald, we conjecture that also
these more general Kostka functions are polynomials in $q$ and
$t^{1/2}$ with non\_negative integers as coefficients. If $q=0$, then
our Kostka functions are Kazhdan\_Lusztig polynomials of a special
type. Therefore, our positivity conjecture combines Macdonald
positivity and Kazhdan\_Lusztig positivity and hints towards a
connection between Macdonald and Kazhdan\_Lusztig theory.

}

\beginsection Intro. Introduction

In \cite{Mac0}, Macdonald introduced a new class of symmetric
functions $J_\mu(z;q,t)$, parameterized by a partition $\mu$ and
depending on two parameters $q$ and $t$, which generalizes both
Hall\_Littlewood and Jack polynomials. In the same paper, he
introduced the two variable Kostka functions
$$68
K_{\lambda\mu}(q,t):=\<s_\lambda,J_\mu\>_{\rm HL}
$$
where $s_\lambda$ is the Schur functions for the partition
$\lambda$ and where $\<\cdot,\cdot\>_{\rm HL}$ denotes the scalar
product rendering the Hall\_Littlewood polynomials orthogonal. Based
on computational evidence and some special cases, Macdonald
conjectured that $K_{\lambda\mu}(q,t)$ is a polynomial in $q$, $t$ with
non\_negative integral coefficients. Even polynomiality was open for a
while and was almost simultaneously proved in \cite{P1}, \cite{P2},
\cite{P3}, \cite{Kostka}, \cite{P4}. Haiman finally
proved positivity in \cite{Hai}.

To prove polynomiality, the author used in \cite{Kostka} a more
general theory, that of non\_symmetric Macdonald polynomials which has
been developed mainly by Cherednik. For that reason, it is tempting to
look for Kostka functions associated to non\_symmetric Macdonald
polynomials and prove their positivity first. In this paper, we
introduce functions $K_{\lambda\mu}(q,t)$ where $\lambda$ and $\mu$
are now allowed to be compositions, i.e., finite unordered sequences
of positive integers and which coincide with Macdonald's when
$\lambda$ and $\mu$ are partitions. This definition links two
theories, Macdonald and Kazhdan\_Lusztig, which, even though they
share the same background, namely affine Hecke algebras, have been
unrelated so far.

The starting point of our theory was Lusztig's observation,
\cite{Lu1}, \cite{Lu2}, that certain Kazhdan\_Lusztig basis elements
can be identified with the Schur function. Therefore, the
idea is roughly to replace in \cite{E68} the symmetric Macdonald
polynomial by a non\_symmetric one and the Schur function
$s_\lambda$ by a Kazhdan\_Lusztig basis element.

More precisely, we consider the standard parabolic module $\cM$ of the
(extended) affine Hecke algebra of type $A_{n-1}$. This module can be
identified with a polynomial ring and has a basis consisting of the
non\_symmetric Macdonald polynomials $\cE_\mu$. On the other hand,
$\cM$ has also the canonical, or Kazhdan\_Lusztig basis
$\uM^\lambda$. This basis is constructed from the standard basis
$M^\lambda$ by forcing selfduality. Moreover, $\cM$ carries the scalar
product on $\cM$ for which the standard basis is orthonormal. Then we
define $K_{\lambda\mu}^{(n)}:=\<\uM^\lambda,\cE_\mu\>$.

So far, the construction works more or less for any root system but
sample calculations show that $K^{(n)}_{\lambda\mu}$ does not have
positive coefficients, even in type $A_{n-1}$. For this to happen, we
have to stabilize, i.e., let the number $n$ tend to $\infty$. If we
equip $\ZZ[q,t^{1/2},t^{-1/2}]$ with the $t$\_adic topology, then we
show that $K_{\lambda\mu}=\|lim|_{n\pfeil\infty}K^{(n)}_{\lambda\mu}$
exists and is an element of $\ZZ[q,t^{1/2},t^{-1/2}]$. This is the
main (proven) result of this paper. 

We conjecture that $K_{\lambda\mu}(q,t)$ has positive coefficients. This
has been confirmed in a great number of cases by direct
computation. Further evidence is the fact that in case $\lambda$ and
$\mu$ are partitions then our $K_{\lambda\mu}$ coincides with
Macdonald's (proved to be positive by Haiman). This is a consequence
of the aforementioned theorem of Lusztig.

Finally, we show that for $q=0$ the Macdonald polynomials specialize
to the standard basis of $\cM$ (see \cite{E69} and also \cite{Kostka}
Cor.~5.7). This means that $K_{\lambda\mu}(0,t)$ is a Kazhdan\_Lusztig
polynomial, hence positive.

At the end of the paper we present a conjecture which substantially
refines Macdonald's positivity. More precisely, we define ``marked''
Kostka polynomials $K_{\lambda\bmu}(t)$ which do not depend on $q$
anymore. Here $\bmu$ is a composition with some boxes of its Young
diagram are marked. The unmarked polynomials $K_{\lambda\mu}$ can be
obtained in an easy and positive way from the marked ones. Ample
numerical evidence suggests that the $K_{\lambda\bmu}(t)$ have
positive coefficients.

\medskip
\noindent{\bf Acknowledgment:} Most of this work was completed 
during a stay at the University of Strasbourg in Spring 1996. It was finished
during a stay at the University of Freiburg. The author thanks both
institutions for their hospitality.

\beginsection Prelim. The affine Hecke algebra for $GL_n$

Consider the lattice $X:=\ZZ^n$ with standard basis
$e_1,\ldots,e_n$. Then the symmetric group on $n$ letters, $\Wf:=S_n$,
acts on $X$ by permuting the $e_i$:
$$
\pi(\tau_1,\ldots,\tau_n):=(\tau_{\pi^{-1}(1)},\ldots,\tau_{\pi^{-1}(n)}).
$$
Let $W:=\Wf\semidir X$ be the (extended) affine Weyl group. For
$\tau\in X$ let $t_\tau$ be the corresponding element in $W$.

Let $X\v:=\|Hom|(X,\ZZ)\cong\ZZ^n$ with basis
$\epsilon_1,\ldots,\epsilon_n$ (the basis dual to
$e_1,\ldots,e_n$). Let $\Delta_f\subseteq X\v$ be the set of roots for
$\Wf$, i.e.,
$$
\Delta_f:=\{\pm(\epsilon_i-\epsilon_j)\mid 1\le i<j\le n\}
$$
We regard $a\in\ZZ$ as the constant function $a$ on $X$. Then the set
of affine roots $\Delta:=\Delta_f+\ZZ$ consists of affine linear
functions on $X$. Let
$$
\alpha_0:=\epsilon_n-\epsilon_1+1,
\alpha_1:=\epsilon_1-\epsilon_2,\ldots,\alpha_{n-1}:=\epsilon_{n-1}-\epsilon_n.
$$
Then $\Sigma_f:=\{\alpha_1,\ldots,\alpha_{n-1}\}$ is the set of simple
roots of $\Delta_f$ while $\Sigma:=\{\alpha_0,\ldots,\alpha_{n-1}\}$ is the
one for $\Delta$. The corresponding simple reflections are denoted by
$s_i$. Thus, for $1\le i<n$ we have $s_i=(i\,i{+}1)$ while
$s_0(\tau_1,\ldots,\tau_n)=(\tau_n+1,\tau_2,\ldots,\tau_{n-1},\tau_1-1)$.
The simple roots generate the positive roots
$$
\Delta_f^+:=\{\epsilon_i-\epsilon_j\mid 1\le i<j\le n\}\subseteq\Delta_f,
\quad
\Delta^+:=\Delta_f^+\cup(\Delta_f+\ZZ_{>0})\subseteq\Delta.
$$
The dominant Weyl chamber is $X_+:=\{\tau\in X\mid \tau_1\ge\ldots\ge\tau_n\}$.

Every element $w\in W$ acts on $\Delta$ by
$w\alpha(\tau)=\alpha(w^{-1}\tau)$. We define its length as
$$
\ell(w):=\#\{\alpha\in\Delta^+\mid w\alpha\in-\Delta^+\}.
$$
For $w\in\Wf$ and $\tau\in X$ we have the useful formula
$$2
\ell(t_\tau w)=\sum_{1\le i<j\le n\atop w^{-1}(i)< w^{-1}(j)}
\big|\,\tau_i-\tau_j\,\big|+
\sum_{1\le i<j\le n\atop w^{-1}(i)>w^{-1}(j)}
\big|\,\tau_i-\tau_j-1\,\big|
$$
This means, in particular, that
$$
\oq:=t_{-e_n}s_{n-1}\ldots s_1=s_{n-1}\ldots s_1t_{-e_1},
$$
acting on $X$ like
$$22
\oq(\tau_1,\ldots,\tau_n)=(\tau_2,\ldots,\tau_n,\tau_1-1),
$$
has length zero. In fact, $\oq$ generates $\Omega:=\{w\in
W\mid\ell(w)=0\}\cong\ZZ$. If $W^a\subseteq W$ is the subgroup
generated by $s_1,\ldots,s_n$, then $W^a=S_n\semidir Q$ with
$Q:=\{\tau\in X\mid \sum_i\tau_i=0\}$ and $W=\Omega\semidir W^a$. The
action of $\Omega$ on $W^a$ is given by
$$
\oq s_i\oq^{-1}=s_{i-1},\quad i=1,\ldots,n-1,\quad
\oq s_0\oq^{-1}=s_{n-1}.
$$

Let $\cL:=\ZZ[v,v^{-1}]$. We often use also the notation $t:=v^2$. The
(extended) affine Hecke algebra $\cH$ is the $\cL$\_algebra generated
by elements $H_0,\ldots,H_{n-1}$, $\omega$ with relations
$$4
H_iH_j=H_jH_i\quad{\rm for}\quad 0\le i<j\le n-1\hbox{ with }1<|i-j|<n-1
$$
$$
H_iH_{i+1}H_i=H_{i+1}H_iH_{i+1} \quad{\rm for}\quad i=0,\ldots,n-1
\quad \hbox{with }H_n:=H_0
$$
and
$$
\oq H_i\oq^{-1}=H_{i-1}\quad i=1,\ldots,n-1,
\quad\oq H_0\oq^{-1}=H_{n-1}
$$
$$1
(H_i+v)(H_i-v^{-1})=0\quad i=0,\ldots,n-1.
$$
If $w=s_{i_1}\ldots s_{i_r}\omega^k\in W$ is a reduced expression, then
one puts $H_w:=H_{i_1}\ldots H_{i_r}\omega^k$. These elements form an
$\cL$\_basis of $\cH$. Moreover
$$5
H_xH_y=H_{xy}\quad{\rm whenever}\ \ell(x)+\ell(y)=\ell(xy)
$$
The subalgebra spanned by $H_w$, $w\in W_f$ is denoted by $\cH_f$.

\beginsection parabolic. The parabolic module

There is a unique $\cL$\_linear homomorphism $\cH_f\pfeil \cL$ with
$H_1,\ldots,H_{n-1}\mapsto v^{-1}$. More generally, $H_w\mapsto
v^{-\ell(w)}$, $w\in W_f$, . This way, $\cL$ becomes a $\cH_f$\_module
denoted by $\cL(v^{-1})$. Consider the induced module
$$
\cM:=\cH\otimes_{\cH_f}\cL(v^{-1}).
$$
Every coset in $W/W_f$ is represented by a unique element $t_\tau$,
$\tau\in X$. This implies that the elements
$H_{t_\tau}\otimes1\in\cM$, $\tau\in X$ form a $\cL$\_basis of
$\cM$. It is convenient to modify this basis slightly. For $\tau\in X$
let $m_\tau$ be the unique shortest element of the coset $t_\tau
W_f$. Using the length formula \cite{E2} one can check that
$m_\tau=t_\tau w_\tau^{-1}$ where
$w_\tau$ is the shortest permutation such that
$w_\tau(\tau)\in -X_+$. A useful formula for $w_\tau$ is
$$52
w_\tau(i)=\#\{j=1,\ldots, i\mid\tau_j\le\tau_i\}+
\#\{j=i+1,\ldots, n\mid\tau_j<\tau_i\}.
$$
The elements
$$40
M_\tau:=m_\tau\otimes1=v^{-\ell(w_\tau)}(t_\tau\otimes1)
$$
form the {\it standard basis\/} of $\cM$.
The action of the generators of $\cH$ in terms
of the standard basis is then given by (see \cite{Soe} \S3)
$$25
(H_i+v)(M_\tau)=\cases{
M_{s_i(\tau)}+vM_\tau&if $\tau_i>\tau_{i+1}$\cr
(v+v^{-1})M_\tau&if $\tau_i=\tau_{i+1}$\cr
M_{s_i(\tau)}+v^{-1}M_\tau&if $\tau_i<\tau_{i+1}$\cr},\qquad i=1,\ldots,n-1
$$
$$26
\oq(M_\tau)=M_{\oq(\tau)}.
$$

The Bruhat order on $W$ induces an order relation on $X$ by defining
$$
\tau\le\eta\Longleftrightarrow m_\tau\le m_\eta.
$$
It has the properties
$$41
s_\alpha(\tau)\ge\tau\Longleftrightarrow\alpha(\tau)\ge0\quad\hbox{for
all}\ \alpha\in\Delta^+
$$
$$53
\tau\le\eta\Longleftrightarrow\omega(\tau)\le\omega(\eta)
$$
$$42
\tau\le\eta
\Longleftrightarrow \|min|\{\tau,s_\alpha(\tau)\}\le s_\alpha(\eta)
\quad
\hbox{for all $\alpha\in\Sigma$ with $\alpha(\eta)\le0$}.
$$
$$54
\tau\le0\Longleftrightarrow\tau=0.
$$
Observe that these properties allow to compute the Bruhat
order algorithmically. In fact, with \cite{E42} one can ``move''
$\eta$ into the fundamental alcove. Then, using \cite{E53}, one
reduces to $\eta=0$. Then one concludes with \cite{E54}.

In general it is not true that $\tau\le\eta$ implies $w\tau\le w\eta$
but there is an important special case when this holds:

\Lemma String. Let $\tau,\eta\in X$, $w\in W_f$ and assume that
$\eta-\tau\in\ZZ\alpha^\vee$ for some $\alpha\in\Delta^+_f$. Assume
moreover $w\alpha>0$. Then $\tau\le\eta$ if and only if $w\tau\le w\eta$.

\Proof: Let $\tau:=\eta-k\alpha^\vee$ and $N:=\alpha(\eta)$. Then
$\tau\le\eta$ if and only if
$$44
k=\cases{0,\ldots,N-1&for $N>0$\cr N,\ldots,0&for $N\le0$\cr}
$$
See, e.g., \cite{Hecke} Lemma 4.1. The result follows since
$w\alpha(w\eta)=\alpha(\eta)$.\qed

\beginsection geometric. The Bernstein presentation

For $\tau\in X_+$ let $X^\tau:=H_{t_\tau}$. If $\tau,\eta\in X_+$, then
\cite{E2} and \cite{E5} imply
$$
X^\tau X^\eta=X^{\tau+\eta}=X^\eta X^\tau.
$$
Hence we can extend the definition for $X^\tau$ to all $\tau\in X$ by
$$
X^\tau:=X^{\tau'}(X^{\tau''})^{-1}\quad\hbox{where }\tau',\tau''\in
X_+\hbox{ with }\tau=\tau'-\tau''.
$$
Since $t_{e_1}=\oq^{-1} s_{n-1}\ldots s_1$ is a reduced expression
this means concretely in our situation $X^\tau=X_1^{\tau_1}\ldots
X_n^{\tau_n}$ with
$$8
X_i:=X^{e_i}=H_{i-1}^{-1}\ldots H_1^{-1}\oq^{-1} H_{n-1}\ldots H_i
$$
This way we get a homomorphism
$$
\Phi:\cL[X]:=\cL[x_1^{\pm1},\ldots,x_n^{\pm1}]\into\cH:x_i\mapsto X_i.
$$
The relations between the $X^\tau$ and the ``finite'' $H_i$
(i.e. $i>0$) have been determined by Bernstein (see \cite{Lu3}~0.3):
$$
H_i\Phi(\xi)-\Phi(s_i(\xi))H_i=(v^{-1}-v)\Phi(x_i{\xi-s_i(\xi)\over
x_i-x_{i+1}}), \qquad i=1,\dots,n-1.
$$
The homomorphism $\Phi$ also identifies $\cL[X]$ with $\cM$ and we get:
$$39
\Psi:\cL[X]\pf\sim\cM:\xi\mapsto\Phi(\xi)\otimes1=\Phi(\xi)(M_0).
$$
By transport of structure, we get an action of $\cH$ on
$\cL[X]$. Concretely, the generators
$H_1,\ldots,H_{n-1},X_1,\ldots,X_n$ act as
$$43
H_i(\xi)=v^{-1}s_i(\xi)+(v^{-1}-v)x_i{\xi-s_i(\xi)\over
x_i-x_{i+1}}, \qquad i=1,\dots,n-1;
$$
$$
X_i(\xi)=x_i\xi, \qquad i=1,\dots,n.
$$
The action of $\oq$ and $H_0$ is more complicated and is deduced
from the relations
$$6
\oq^{-1}=H_1\ldots H_{n-1}X_n=X_1H_1^{-1}\ldots H_{n-1}^{-1}.
$$
$$7
H_0=\oq^{-1}H_{n-1}\oq=\oq H_1\oq^{-1}=
X_1X_n^{-1}H_{(1n)}^{-1}\quad\hbox{with } (1n):=s_1\ldots s_{n-1}\ldots s_1.
$$

\Lemma triang. For every $w\in W_f$ and $\tau\in X$ holds
$$45
H_w(x^\tau)\in v^{\ell(w)-2k}x^{w(\tau)}+
\sum_{\eta<w(\tau)}\cL x^\eta.
$$
where
$$46
k:=\#\{\alpha\in\Delta^+_f\mid w\alpha<0,\ \alpha(\tau)\ge0\}.
$$

\Proof: If $w=1$, the statement is trivial. Otherwise write $w=vs$
with $s$ a simple reflection and $\ell(v)<\ell(w)$. This means
$v\alpha_s>0$. Let $x^\mu$ be a monomial occurring in $H_s(x^\tau)$ and
$x^\eta$ a monomial occurring in $H_v(x^\mu)$. By the explicit formula
\cite{E43} and \cite{E44} we get $\mu\le s(\tau)$ and
$s(\tau)-\mu\in\ZZ\alpha_s^\vee$. Then \cite{String} implies
$v(\mu)\le w(\tau)$. By induction, we have $\eta\le v(\mu)$, hence
$\eta\le w(\tau)$. Finally, the coefficient of $x^{s\tau}$ in
$H_s(x^\tau)$ is $v^{-1}$ or $v$ according to $\alpha_s(\tau)\ge0$ or
$\alpha_s(\tau)<0$, respectively. On the other hand, $k=k_w(\tau)$
satisfies the recursion
$$
k_w(\tau)=k_v(s\tau)+\cases{1&if $\alpha_s(\tau)\ge0$\cr0&otherwise\cr}
$$
which implies the claim on the leading coefficient.
\qed

\beginsection KazhdanLusztig. The Kazhdan-Lusztig basis

The defining relations \cite{E4}--\cite{E1} of $\cH$ imply that it
admits a unique ring automorphism $d$ with
$$
d(v)=v^{-1},\quad d(\oq)=\oq,\quad d(H_i)=H_i^{-1}\quad i=1,\ldots,n.
$$
More generally we have $d(H_w)=H_{w^{-1}}^{-1}$ for any $w\in W$. Now
put $\cH_{++}:=\sum_{x\in W}v\ZZ[v]H_x$. Then we have the following
fundamental result of Kazhdan\_Lusztig (\cite{KL}, see also
\cite{Soe}):

\Theorem. For every $w\in W$ there is a unique $\uH_w\in\cH$ with
$d(\uH_w)=\uH_w$ and $\uH_w\in H_w+\cH_{++}$. This element is
triangular with respect to the Bruhat order, i.e., $\uH_w\in\sum_{v\le
w}\cL H_v$. Moreover, the collection of $\uH_w$, $w\in W$ forms an
$\cL$\_basis of $\cH$.

A similar construction works for $\cM$. Since the homomorphism
$\cH_f\pfeil\cL$ defining $\cL(v^{-1})$ commutes with $d$ we may
define an involution of $\cM$, also denoted by $d$, by
$$
d(\xi\otimes a(v)):=d(\xi)\otimes a(v^{-1}).
$$
Again we form $\cM_{++}:=\sum_{\tau\in X}v\ZZ[v]M_\tau$ and obtain
(see \cite{Soe}):

\Theorem. For every $\tau\in X$ there is a unique $\uM_\tau\in\cM$
with $d(\uM_\tau)=\uM_\tau$ and $\uM_\tau\in M_\tau+\cM_{++}$. This
element is triangular with respect to the Bruhat order, i.e.,
$\uM_\tau\in\sum_{\eta\le \tau}\cL M_\eta$. Moreover, the collection
of $\uM_\tau$, $\tau\in X$ forms an $\cL$\_basis of $\cM$.

\Remark: Observe that the triangularity property of $\uM_\tau$ implies
easily that also $d$ is triangular, i.e.,
$d(M_\tau)\in\sum_{\eta\le\tau}\cL M_\eta$.  \medskip Using the
bijection $\Psi:\cL[X]\pf\sim\cM$ the involution $d$ may be
transported to $\cL[X]$, which we also denote by $d$. Explicitly we get
(\cite{Hecke} Lemma~3.4):

\Theorem involX. Let $p\mapsto\pq$ be the ring involution of $\cL[X]$ with
$\vq=v^{-1}$ and $\overline{x_i}=x_i$. Let $w_0\in W_f$ be the longest
element. Then for any $p\in\cL[X]$ holds
$$56
d(\Phi(p))=H_{w_0}\Phi(w_0\pq)H_{w_0}^{-1}.
$$
In particular,
$$23
d(p)=v^{\ell(w_0)}H_{w_0}(w_0\pq)).
$$

\noindent We are going to need only two properties of the
Kazhdan\_Lusztig basis.

\Lemma KLsymmetry. Let $\tau\in X$ and $s\in\Sigma$ with
$s(\tau)\le\tau$. Then $H_s(\uM_\tau)=v^{-1}\uM_\tau$.

\noindent For a proof see, e.g., \cite{Soe} Prop.~3.6. The Lemma
implies in particular that $\Psi^{-1}(\uM_\tau)$ is symmetric whenever
$\tau\in-X_+$. In fact, it can be computed explicitly:

\Theorem Lusztig. For $\lambda\in X_+$ let $s_\lambda$ be the
corresponding Schur polynomial. Then
$$
\uM_{-\lambda}=\Psi(s_\lambda(x_1^{-1},\ldots,x_n^{-1}))
$$

\noindent This result is due to Lusztig, first proved for $A_{n-1}$ in
\cite{Lu1} and then for arbitrary root systems in \cite{Lu2}. For an
alternate proof see \cite{Hecke}. It is the key to our approach to
Kostka polynomials.

\beginsection Macdonald. Macdonald polynomials

As mentioned, the action \cite{E6} of $\omega$ on $\cL[X]$ is quite
complicated. Cherednik had the idea (see e.g. \cite{Ch}) to replace
$\omega$ by
$$
\tilde\oq(f)(x_1,\ldots,x_n):=f(qx_n,x_1,\ldots,x_{n-1})
$$
where $q$ is an additional parameter. This formula is motivated by the
affine linear action \cite{E22} of $\oq$. Also the action of $H_0$ becomes
easy this way:
$$
\tilde H_0:=\tilde\oq H_1\tilde\oq^{-1}=
v^{-1}\tilde s_0+(v^{-1}-v)x_n{1-\tilde s_0\over x_n-q^{-1}x_1}
$$
with
$$
\tilde s_0(p):=p(qx_n,x_2,\ldots,x_{n-1},q^{-1}x_1).
$$
One checks that $\tilde H_0,H_1,\ldots,H_{n-1},\tilde\oq$ satisfy the
relations \cite{E4}--\cite{E1} and therefore generate another copy
$\tilde\cH$ of $\cH$. In particular, $\tilde\cH$ will contain a copy
of $\cL[X]$ which we choose to be generated by the elements
$$10
\xi_i:=v^{1-n}H_{i-1}\ldots H_1\tilde\oq^{-1} H_{n-1}^{-1}\ldots H_i^{-1},
\quad i=1,\ldots,n
$$
Note that this definition is ``dual'' to \cite{E8} and also has the
factor $v^{1-n}$. The reason for this is to get later the stability
property \cite{E65}. The main feature of $\tilde\cH$ is that it acts
locally finitely on $\cL_q[X]$ where $\cL_q:=\cL[q,q^{-1}]$. More
precisely,

\Lemma. For $i=1,\ldots,n$ and $\tau\in X$ holds
$$
\xi_i(x^\tau)\in q^{-\tau_i}t^{1-w_\tau(i)}\,
x^\tau+\sum_{\mu<\tau}\cL_qx^\mu.
$$

\Proof: First, we show that the $\xi_i$ are triangular with respect to
the Bruhat order. It suffices to do this for
$\Xi_i:=\xi_i\xi_{i+1}\ldots\xi_n$ with $i=1,\ldots,n$. Formula
\cite{E10} implies $\Xi_i=v^{-N}H_{w_i}\tilde\oq^{i-n-1}$ where
$w_i(\tau)=(\tau_{i+1},\ldots,\tau_n,\tau_1,\ldots,\tau_i)$ and
$N=(n-i+1)(n-1)$. Thus,
$$
\Xi_i(x^\tau)=v^{-N}H_{w_i}\tilde\oq^{i-n-1}(x^\tau)=
q^{-\tau_i-\ldots-\tau_n}v^{-N}H_{w_i}x^{w_i^{-1}(\tau)}\in
\sum_{\mu\le\tau}\cL_q x^\mu
$$
by \cite{triang}. The formula for the leading coefficient follows easily
from \cite{E45}, \cite{E46}, and \cite{E52}.\qed

\noindent
Since the Cherednik operators $\xi_1,\ldots,\xi_n$ commute and are
triangular with distinct diagonal terms they have a common eigenbasis,
the {\it non-symmetric Macdonald polynomials} $\cE_\lambda$.

\Corollary. For every $\lambda\in X$ there is $\cE_\lambda\in\cL_q[X]$,
unique up to a scalar, with
$$13
\xi_i(\cE_\lambda)=q^{\lambda_i}t^{1-w_{-\lambda}(i)}\cE_\lambda,\quad
i=1,\ldots,n
$$
Moreover, $\cE_\lambda$ is triangular with respect to the Bruhat
order:
$$47
\cE_\lambda\in\sum_{\mu\le-\lambda}\cL_q x^\mu.
$$

\Remark: Usually (see, e.g., \cite{Mac2} (2.7.5)), the triangularity
of $\cE_\lambda$ is expressed with respect to an order which is finer
than the Bruhat order.
\medskip

We are normalizing $\cE_\lambda$ in the following way. As usual, we
represent $\lambda$ by its {\it diagram}, i.e., the set of pairs
$(i,j)\in\ZZ^2$ (called {\it boxes}) with $1\le j\le\lambda_i$. To a
box $s=(i,j)\in\lambda$ we associate its {\it arm-length}
$$
a_\lambda(s):=\lambda_i-j
$$
and its {\it leg\_length}
$$
l_\lambda(s):=\#\{k<i\mid j\le\lambda_k+1\le\lambda_i\}+
\#\{k>i\mid j\le\lambda_k\le\lambda_i\}.
$$
Now we demand that the coefficient of $x^{-\lambda}$ in $\cE_\lambda$ is
$$55
\prod_{s\in\lambda}\left(1-q^{a_\lambda(s)+1}t^{l_\lambda(s)+1}\right).
$$
One can show, \cite{Kostka} Cor.~5.2, that with this normalization the
coefficients of $\cE_\lambda$ are polynomials in $q$ and $t$.

\beginsection polynomial. The polynomial part of $\cM$

Subsequently, we are only interested in the ``polynomial'' part of
$\cM$. The reason for this is its stability properties as
$n\pfeil\infty$. Let us first introduce the polynomial part $\Hp$
of $\cH$, namely the subalgebra generated by $\cH_f$, and $\oq$ (but
{\it not} $\oq^{-1}$). Let $\Lambda:=\NN^n\subseteq X$ and consider
the submonoid $\Wp:=W_f\semidir(-\Lambda)$ of $W$. Then we have:

\Theorem. The set $\{H_w\mid w\in\Wp\}$ is an $\cL$\_basis of
$\Hp$. Moreover, $Z_i:=X_i^{-1}\in\Hp$ for all $i$ and
$$38
\cL[Z_1,\ldots,Z_n]\otimes_\cL\cH_f\pfeil\Hp: p(Z)\otimes u\mapsto p(Z)u
$$
is bijective. Furthermore, $\Hp$ has the following presentations:
\Item{1}It is generated by the subalgebras $\cH_f$ and $\cL[Z_1,\ldots,Z_n]$
with relations
$$28
\displaystyle{H_ip-(s_ip)H_i=
(v-v^{-1})Z_{i+1}{p-(s_ip)\over Z_i-Z_{i+1}}}
\ {\rm for}\ i=1,\ldots,n-1, p\in\cL[Z_1,\ldots,Z_n].
$$
\Item{2}It is generated by the subalgebra $\cH_f$ and $Z_1,\ldots,Z_n$
with relations
$$0
\eqalignno{%
29&&H_iZ_iH_i=Z_{i+1}\ {\rm for}\  i=1,\ldots,n-1\cr
30&&H_jZ_i=Z_iH_j\ {\rm for}\  i=1,\ldots,n, j=1,\ldots,n-1, i-j\ne0,1\cr
31&&Z_iZ_j=Z_jZ_i\ {\rm for}\  i,j=1,\ldots,n\cr}
$$
\Item{3}It is generated by the subalgebra $\cH_f$ and $Z_1$ with relations
$$0
\eqalignno{%
32&&H_1Z_1H_1Z_1=Z_1H_1Z_1H_1\cr
33&&H_iZ_1=Z_1H_i\ {\rm for}\ i=2,\ldots,n\cr}
$$
\Item{4}It is generated by the subalgebra $\cH_f$ and $\oq$ with relations
$$0
\eqalignno{%
34&&H_i\oq=\oq H_{i+1}\ {\rm for}\ i=1,\ldots,n-2\cr
35&&H_{n-1}\oq^2=\oq^2H_1\cr}
$$\Par

\Proof: Formula \cite{E8} implies $Z_i\in\Hp$. Let $w_0\in W_f$ be the
longest element. Then, for
$\lambda\in\Lambda\cap X_+$
$$
\Hp\ni Z^\lambda=(X^\lambda)^{-1}=H_{t_\lambda}^{-1}=d(H_{t_{-\lambda}})=
H_{w_0}H_{t_{-w_0(\lambda)}}H_{w_0}^{-1}
$$
(by \cite{E56}). This implies $H_{t_{-\lambda}}\in\Hp$ for all
$\lambda\in\Lambda\cap X_+$. Since $H_{sw}=H_s^{\pm1}H_w$ and
$H_{ws}=H_wH_s^{\pm1}$ for all simple reflections $s$ we get
$\cH'\subseteq\Hp$ where $\cH'$ is the $\cL$\_submodule spanned by all
$H_w\in\Hp$, $w\in\Wp$. Conversely, formula \cite{E5} implies
$\omega\cH'\subseteq\cH'$. Hence $\Hp\subseteq\cH'$.

\cite{I1} Let now $\cH'$ be the algebra generated by $\cL[Z]$ and
$\cH_f$ subject to relation \cite{E28}. Then there are natural maps
$$
\matrix{\cL[Z]\otimes_\cL\cH_f&\Pf{\phi_1}&\cH'&\Pf{\phi_2}&\Hp\cr
\downarrow\Rechts{\phi_3}&&&&\downarrow\Rechts{\phi_4}\cr
\cL[Z,Z^{-1}]\otimes_\cL\cH_f&&\Pf{\phi_5}&&\cH\cr}.
$$
Since $\phi_3$ and $\phi_5$ are injective, also $\phi_1$ is
injective. Relation \cite{E28} implies that $\phi_1$ is also
surjective. Thus, $\phi_2$ is injective. Now $\omega=H_{n-1}\ldots
H_1Z_1$ implies that $\phi_2$ is also surjective. This implies the
bijectivity of \cite{E38}.

\cite{I2} Relation \cite{E31} simply means that $\Hp$ contains
$\cL[Z]$ as a subalgebra. Moreover \cite{E29}, \cite{E30} are
equivalent to \cite{E28} for $p=Z_j$. It is well known (see
\cite{Lu3}~3.6) that that case
implies \cite{E28} for any $p$. Thus, the presentation \cite{I1} and
\cite{I2} are equivalent.

\cite{I3} Here, we are defining $Z_2,\ldots,Z_n$ using formula
\cite{E29}. Then \cite{E32} is nothing else than $Z_2Z_1=Z_1Z_2$. Thus,
\cite{I2} implies \cite{I3}. Conversely, assume \cite{E32}, \cite{E33}
hold. Relations \cite{E29} are true by definition. Relation \cite{E30}
follows easily for $i<j$. For $i=j+2$ we have
$$
\eqalign{H_jZ_i&=H_jH_{j+1}H_jZ_jH_jH_{j+1}=
H_{j+1}H_jH_{j+1}Z_jH_jH_{j+1}=\cr
&=H_{j+1}H_jZ_jH_{j+1}H_jH_{j+1}=
H_{j+1}H_jZ_jH_jH_{j+1}H_j=
Z_iH_j.\cr}
$$
For $i>j+2$ we get \cite{E30} by induction from
$Z_i=H_{i-1}Z_{i-1}H_{i-1}$. Finally, \cite{E31} follows the same way
by induction from \cite{E30} and \cite{E32}.

\cite{I4} First assume the relations in \cite{I4}. Define
$Z_1:=H_1^{-1}\ldots H_{n-1}^{-1}\omega$. Then
$$0
\eqalignno{&H_1Z_1H_1Z_1&=H_2^{-1}\ldots H_{n-1}^{-1}\omega H_2^{-1}\ldots
H_{n-1}^{-1}\omega=\cr
36&&=(H_2^{-1}\ldots H_{n-1}^{-1}
H_1^{-1}\ldots H_{n-2}^{-1}H_{n-1}^{-1})H_{n-1}\omega^2\cr}
$$
$$0
\eqalignno{&Z_1H_1Z_1H_1&=H_1^{-1}\ldots H_{n-1}^{-1}\omega
H_2^{-1}\ldots H_{n-1}^{-1}\omega H_1=\cr
37&&=(H_1^{-1}\ldots H_{n-1}^{-1}
H_1^{-1}\ldots H_{n-2}^{-1})\omega^2H_1\cr}
$$
Both expressions are equal since the terms in parenthesis correspond to reduced
expressions of the same permutation namely $(n,n-1,1,\ldots,n-2)$. Moreover,
for $i\ge2$ we have
$$
\eqalign{
H_i^{-1}Z_1
&=H_1^{-1}\ldots H_{i-2}^{-1}(H_i^{-1}H_{i-1}^{-1}H_i^{-1})H_{i+1}^{-1}
\ldots H_{n-1}^{-1}\omega=\cr
&=H_{i-2}^{-1}(H_{i-1}^{-1}H_i^{-1}H_{i-1}^{-1})H_{i+1}^{-1}
\ldots H_{n-1}^{-1}\omega=Z_1H_i^{-1}.\cr}
$$
Assume now conversely that relations \cite{I2} hold. We define
$\omega:=H_{n-1}\ldots H_1 Z_1$. Then
$$
\omega H_{i+1}=H_{n-1}\ldots H_{i+1}H_iH_{i+1}H_{i-1}\ldots
H_1Z_1=H_i\omega
$$
Finally, \cite{E35} can be deduced from the equality of \cite{E36} and
\cite{E37}.\qed

\Remark: The proof shows that \cite{I2}--\cite{I4} are also equivalent
presentations of a ``braid monoid'' with generators
$H_1^{\pm1},\ldots,H_n^{\pm1}$ and $\omega$. Of course, the braid
relations among the $H_i$ should also hold.
\medskip

Now let $\Mp$ be the $\cL$\_submodule of $\cM$ spanned by all
$M_\lambda$ with $\lambda\in-\Lambda$. 

\Theorem. The map
$$
\Hp\otimes_{\cH_f}\cL(v^{-1})\pfeil\cM
$$
is injective with image $\Mp$. In particular, $\Mp$ is an
$\Hp$\_module and
$$
\Psi:\cL[z_1,\ldots,z_n]\pf\sim\Mp
$$
is an isomorphism where $z_i:=x_i^{-1}$.

\Proof: Consider the following commutative diagram:
$$
\matrix{\cL[z]&\pf{\phi_1}&\Hp\otimes_{\cH_f}\cL(v^{-1})\cr
\untenPf\Rechts{\phi_2}&&\untenPf\Rechts{\phi_3}\cr
\cL[z,z^{-1}]&\pf{\phi_4}&\cM\cr}
$$
Then $\phi_1$, $\phi_4$ are bijective by \cite{E38}, \cite{E39},
respectively, while $\phi_2$ is obviously
injective. Hence, $\phi_3$ is injective. Formulas \cite{E25},
\cite{E26} show that $\Mp$ is an $\Hp$\_module. Hence
$\|Im|\phi_3\subseteq\Mp$. The converse inclusion follows from
\cite{E40}.\qed

Observe that the operators $H_i$, $i=1,\ldots,n-1$ and $\tilde\oq$ take
the following form in the coordinates $z_i=x_i^{-1}$:
$$18
H_i=v^{-1}s_i+(v-v^{-1})z_{i+1}{1-s_i\over z_i-z_{i+1}}
$$
$$66
\tilde\oq(f)(z_1,\ldots,z_n)=f(q^{-1}z_n,z_1,\ldots,z_{n-1})
$$
Moreover, to simplify notation, we write from now on for $\lambda\in\Lambda$
$$
M^\lambda:=M_{-\lambda},\ \uM^\lambda:=\uM_{-\lambda},\
w^\lambda:=w_{-\lambda}.
$$
Observe that $w^\lambda$ is the shortest permutation such that
$w^\lambda(\lambda)$ is a partition. We also modify~$\omega$:
$$
\omega^*(\tau):=-\omega(-\tau)=(\tau_2,\ldots,\tau_n,\tau_1+1)\quad
\hbox{hence}\quad \omega(M^\lambda)=M^{\omega^*(\lambda)}
$$
Finally, the modified Bruhat order is
$$
\lambda\preceq\mu\Longleftrightarrow -\lambda\le-\mu.
$$

\Lemma closed. Let $\lambda\in X$ and $\mu\in\Lambda$. Then
$\lambda\preceq\mu$ implies $\lambda\in\Lambda$.

\Proof: Let $\lambda,\mu$ be a counterexample with $N:=|\mu|=\sum_i\mu_i$
minimal. Let $w\in W_f$ be minimal with $\mu':=w\mu\in -X_+$, i.e.,
$\mu'_1\le\ldots\le\mu_n'$. By \cite{E41}, \cite{E42} there is $w'\in
W_f$, $w'\le w$ with $\lambda':=w'\lambda\preceq\mu'$. Clearly also
$\lambda'\not\in\Lambda$ and $\mu\ne0$. Hence $\mu'_n>0$ which implies
that
$\mu'':=(\omega^*)^{-1}(\mu')=(\mu'_n-1,\mu'_1,\ldots,\mu'_{n-1})\in\Lambda$
while $\lambda'':=(\omega^*)^{-1}(\lambda')\not\in\Lambda$. Since
$\lambda''\preceq\mu''$ and $\sum_i\mu''_i=N-1$ we get a contradiction to
the minimality of $N$.\qed

\Corollary. The subset $\Mp$ is stable under the involution
$d$. Moreover, the Kazhdan\_Lusztig elements $\uM^\lambda$ with
$\lambda\in\Lambda$ form an $\cL$\_basis of $\Mp$.

\Proof:  This follows from
\cite{closed} and the triangularity of the involution $d$ and the
Kazhdan\_Lusztig elements.\qed

\Corollary. For $\lambda\in\Lambda$, the Macdonald polynomial
$\cE_\lambda$ is in $\Mp_q:=\Mp\otimes_\cL\cL_q$.

\Proof: Follows immediately from \cite{closed} and the triangularity
property \cite{E47}.\qed

\beginsection recursion. Recursion formulas for Macdonald polynomials

In this section, we describe the recursion formulas from \cite{Kostka}%
\footnote{For the convenience of the reader we include following
conversion table between notations:\medskip
{\offinterlineskip
\halign{\vrule\ #\hfill\ \vrule width.8pt&&\ $#$\hfill\ \vrule\cr
\noalign{\hrule}
\vrule height 10pt width 0pt depth 4pt\cite{Kostka}&t&H_i&\Hq_i&
\Delta&\xi_i&\cE_\lambda&\Phi&A_m&\Aq_m&\Phi'&A'_m&\Aq'_m\cr
\noalign{\hrule}
\vrule height 12pt width 0pt depth 4ptThis paper&t=v^2&vH_i^{-1}&vH_i&\tilde\omega&
\xi_i&\cE_\lambda&Z_n\tilde\omega&v^{n-m}\overline{\tilde\Phi}_m&
v^{n-m}\Phi_m&v^{1-n}\omega&v^{1-m}\overline\Phi_m&v^{1-m}\Phi_m\cr
\noalign{\hrule}}}}
which produce exactly the $\cE_\lambda$ with $\lambda\in\Lambda$. For
$m=1,\ldots,n$ we define the operators
$$11
\tilde\Phi_m:=H_m...H_{n-1}Z_n\tilde\oq,
$$
$$12
\overline{\tilde\Phi}_m:=H_m^{-1}...H_{n-1}^{-1}Z_n\tilde\oq.
$$
Recall that the length $l(\lambda)$ of $\lambda\in \Lambda=\NN^n$ is
the maximal $m\ge0$ with $\lambda_m\ne0$ (so $\lambda=0$ if and only
if $l(\lambda)=0$).

\Theorem. {\rm (\cite{Kostka} Thm.~5.1)} For $\lambda\in \Lambda$ let
$m:=l(\lambda)$. Then
$$9
\cE_\lambda=
q^{\lambda_m-1}v^{n-m}\left(\tilde\Phi_m-
q^{\lambda_m}v^{2a}\overline{\tilde\Phi}_m\right)(\cE_{\lambda^*})
\quad\hbox{where}
$$
$$0
\eqalignno{&\lambda^*&
:=(\lambda_m-1,\lambda_1,\ldots,\lambda_{m-1},0,\ldots,0)\cr
59&a&:=1+\#\{i=1,\ldots,m\mid\lambda_i<\lambda_m\}.\cr}
$$

\medskip\noindent Clearly, starting from $\cE_0=1$, this formula allows to
compute $\cE_\lambda$ for all $\lambda\in \Lambda$ in a unique way.

Following \cite{Kostka}, we are going to rewrite the recursion
\cite{E9}. Equations \cite{E8} and \cite{E10} imply
$$
Z_n\tilde\oq=v^{1-n}\oq\xi_1^{-1}.
$$
Inserting this into \cite{E11}, \cite{E12} and observing that
$\cE_{\lambda^*}$ is an eigenvector for $\xi_1^{-1}$ (see \cite{E13}) we
obtain
$$19
\cE_\lambda=
v^{m+1-2a}\left(\Phi_m-
q^{\lambda_m}v^{2a}\overline\Phi_m\right)(\cE_{\lambda^*})
$$
with the new operators
$$
\Phi_m:=H_m...H_{n-1}\oq=H_{c_m},
$$
$$
\overline\Phi_m:=H_m^{-1}...H_{n-1}^{-1}\oq=H_{c_m^{-1}}^{-1}=d(H_{c_m})
$$
where
$$
c_m:=s_{m-1}\ldots s_1t_{-e_1}:\lambda\mapsto
(\lambda_2,\ldots,\lambda_m,\lambda_1-1,\lambda_{m+1},\ldots,\lambda_n)
$$
For the renormalized Macdonald polynomial
$$14
\tilde\cE_\lambda=v^{\ell(w^\lambda)}\cE_\lambda
$$
we obtain the simple formula
$$21
\tilde\cE_\lambda=
\left(\Phi_m-q^{\lambda_m}t^a\overline\Phi_m\right)(\tilde\cE_{\lambda^*}).
$$
The big advantage of \cite{E19} and \cite{E21} over \cite{E9} is that the
parameter $q$ is not involved in the operators $\Phi_m$ and $\overline\Phi_m$.

The first application of the recursion formulas is the following
integrality result from \cite{Kostka}:

\Theorem. For every $\lambda\in \Lambda$ holds
$\cE_\lambda\in\ZZ[t,q][z_1,\ldots,z_n]$.

\Proof: The definitions \cite{E11},\cite{E12} and formula \cite{E18}
show that the operators $v^{n-m}\tilde\Phi_m$ and
$v^{n-m}\overline{\tilde\Phi}_m$ preserve the ring
$\ZZ[t,q,q^{-1}][z_1,\ldots,z_n]$. Thus, formula \cite{E9} implies
that $\cE_\lambda$ is in this ring. On the other hand, \cite{E19}
implies clearly the non\_occurrence of negative powers of $q$.\qed

\noindent
When we express $\cE_\lambda$ in terms of the standard basis then we get

\Corollary. $\Psi(\cE_\lambda)\in\sum\limits_{\mu\in\Lambda}
\ZZ[v^2,q]v^{-\ell(w^\mu)}M^\mu$.

\Proof: \cite{Hecke} Lemma~4.2 implies that the transition
matrix between monomials $z^\tau$ and the elements
$v^{-\ell(w^\mu)}\Psi^{-1}(M^\mu)$ is unitriangular with coefficients
in $\ZZ[t]$.\qed

\noindent According to this Corollary, the coefficients of
$\Psi(\cE_\lambda)$ might contain arbitrary large negative powers of
$v$. Computational evidence leads to:

\Conjecture Vermutung1.
$\Psi(\tilde\cE_\lambda)\in\sum_{\mu\in\Lambda}\ZZ[v,q]M^\mu$.

Another consequence of the recursion formula is that $\cE_\lambda$ is
almost selfdual. More precisely, using the isomorphism
$\Psi:\cL[Z_i]\pf\sim\Mp$ we can transport the involution $d$ to
$\cL[Z_i]$ (\cite{involX}). Then we extend it to $\cL_q[Z_i]$ by
defining $d(q):=q^{-1}$.

\Theorem. For every $\lambda\in\Lambda$ holds
$$
d(\cE_\lambda)=(-1)^{|\lambda|}q^{-A}t^{-B}\cE_\lambda,\quad
d(\tilde\cE_\lambda)=
(-1)^{|\lambda|}q^{-A}t^{-B-\ell(w^\lambda)}\tilde\cE_\lambda
$$
where
$$
A:=\sum_{i\ge1}{\lambda_i+1\choose2}
$$
$$
B:=\sum_{i\ge1}i\lambda_i^+\qquad\hbox{with}\quad \lambda^+
:=w^\lambda(\lambda)\in X_+.
$$

\Proof: The formula for $\tilde\cE_\lambda$ follows immediately from
that for $\cE_\lambda$ and the definition \cite{E14}.

Write $A(\lambda)$ and $B(\lambda)$ for $A$ and $B$, respectively. We
proceed by induction on $|\lambda|$. The assertion is obvious for
$\lambda=0$. Now assume it holds for $\lambda^*$. Clearly, we have
$d\Phi_m=\overline\Phi_md$. Hence
$$
\eqalign{
d(\cE_\lambda)&=v^{-m-1+2a}(\overline\Phi_m-q^{-\lambda_m}v^{-2a}\Phi_m)
((-1)^{|\lambda^*|}q^{-A(\lambda^*)}v^{-2B(\lambda^*)}\cE_{\lambda^*})=
\cr&=
(-1)^{|\lambda^*|+1}q^{-A(\lambda^*)-\lambda_m}v^{-2B(\lambda^*)-m-1}
(-q^{\lambda_m}v^{2a}\overline\Phi_m+\Phi_m)(\cE_{\lambda^*})=
\cr&=
(-1)^{|\lambda^*|+1}q^{-A(\lambda^*)-\lambda_m}v^{-2B(\lambda^*)-2(m+1-a)}
\cE_\lambda
\cr}
$$
Thus we have to show
$$15
|\lambda|=|\lambda^*|+1,
$$
$$16
A(\lambda)=A(\lambda^*)+\lambda_m,
$$
$$17
B(\lambda)=B(\lambda^*)+(m+1-a).
$$
Equation \cite{E15} is obvious, \cite{E16} is easy, and it remains to
prove \cite{E17}. Let $\lambda_k^+$ be the rightmost entry of
$\lambda^+$ which equals $\lambda_m$. Then
$$
k=\#\{i=1,\ldots,m\mid\lambda_i\ge\lambda_m\}=m-(a-1)=m+1-a.
$$
On the other hand $(\lambda^*)^+$ differs from $\lambda^+$ only in its
$k$\_th entry which is $\lambda_m-1$. Hence
$B(\lambda)=B(\lambda^*)+k$ which proves \cite{E17}.\qed

\beginsection Stabilization. Stabilization

Now we want to study Macdonald and Kazhdan\_Lusztig polynomials as
$n\pfeil\infty$. The Hecke algebra studied so far will be denoted by
$\cH_n$. Its parabolic module is $\cM_n$ with its polynomial subset
$\Mp_n$. The element $\omega$ of $\cH_n$ will be denoted $\omega_n$.

\Theorem commrel. Let $\pi_n:\Mp_n\pfeil\Mp_{n-1}$ be the projection with
$$51
\pi(M^\lambda)=\cases{M^{\lambda'}&if $\lambda_n=0$ and where
$\lambda':=(\lambda_1,\ldots,\lambda_{n-1})$\cr
0&if $\lambda_n>0$.\cr}
$$
Then the following commutation relations hold:
$$48
\pi_n H_i=H_i\pi_n\quad i=1,\ldots,n-2
$$
$$50
\pi_n\oq_n=0
$$
$$49
\pi_n H_{n-1}\oq_n=\pi_n H_{n-1}^{-1}\oq_n=\oq_{n-1}\pi
$$
$$24
\pi_n Z_i=\cases{Z_i\pi_n&for $i=1,\ldots,n-1$\cr
0&for $i=n$\cr}
$$

\Proof: Equation \cite{E48} follows immediately from \cite{E25}. For
$\lambda\in\Lambda$ let
$\lambda^*:=\omega_n^*(\lambda)=(\lambda_2,\ldots,\lambda_n,\lambda_1+1)$.
Then \cite{E50} follows from $\omega_n(M^\lambda)=M^{\lambda^*}$. Moreover
$$
\pi_nH_{n-1}\omega_n(M^\lambda)=
\pi_nH_{n-1}M^{\lambda^*}\quad\cases{
=\pi_n M^{s_{n-1}(\lambda^*)}=\omega_{n-1}\pi_n(M^\lambda)
&if $\lambda_n=0$\cr
\in\cL \pi_n M^{s_{n-1}(\lambda^*)}+\cL \pi_n M^{\lambda^*}=0
& if $\lambda_n>0$\cr}
$$
This proves the first part of \cite{E49}. The second part follows
using \cite{E50}. Finally, we get \cite{E24} by using the above and the
explicit expression \cite{E8} for $Z_i=X_i^{-1}$.\qed

\noindent Let
$\pi_n:\cL[z_1,\ldots,z_{n-1},z_n]\pfeil\cL[z_1,\ldots,z_{n-1}]$ be
the obvious projection. Then equation \cite{E24} implies:

\Corollary commdia. The following diagram commutes
$$
\matrix{\cL[z_1,\ldots,z_{n-1},z_n]&\pf{\Psi_n}&\Mp_n\cr
\untenPf\Rechts{\pi_n}&&\untenPf\Rechts{\pi_n}\cr
\cL[z_1,\ldots,z_{n-1}]&\pf{\Psi_{n-1}}&\Mp_{n-1}\cr}
$$

Both $\cL[z_1,\ldots,z_{n-1},z_n]$ and $\Mp_n$ carry a
natural grading, the first by degree, the second by defining
$\|deg|M^\lambda:=|\lambda|=\sum_i\lambda_i$. Moreover, $\Psi_n$ is
degree\_preserving. This follows from the definition \cite{E8} of
$Z_i=X_i^{-1}$ and the fact that $H_i$, $\omega$ is homogeneous of
degree $0$ and $1$, respectively (see \cite{E25}, \cite{E26}). 

\cite{commdia} implies that if we consider the projective limits
$$
\Mp_\infty:=\plim\Mp_n,\quad\cP_\infty:=\plim\cL[z_1,\ldots,z_n]
$$
in the category of graded abelian groups, then we get an isomorphism
$$
\Psi:\cP_\infty\pf\sim\Mp_\infty
$$
More precisely, let $\Lambda:=\NN^{(\infty)}$ be the set of all
sequences of natural numbers almost all of which are zero. Then
$\Mp_\infty$, $\cP_\infty$ is the set of all possibly infinite sums
$\sum_\lambda a_\lambda M^\lambda$, $\sum_\lambda a_\lambda
z^\lambda$, respectively where $\lambda$ runs through a subset of
$\Lambda$ in which $|\lambda|$ remains bounded.

A further consequence of \cite{commrel} is

\Corollary. The space $\Mp_\infty$
carries an action of the operators $H_i$, $Z_i$, $\Phi_i$,
$\overline\Phi_i$ $(i\ge1)$.

Next we need a property of the Bruhat order:

\Lemma stableBruhat. Fix an $i$ with $1\le i\le n+1$. For $\lambda\in
\ZZ^{n+1}$ let $\lambda'\in\ZZ^n$ be obtained from $\lambda$ by
omitting the $i$\_th entry. Let $\lambda,\mu\in\ZZ^n$ with
$\lambda_i=\mu_i$. Then $\lambda\le\mu$ if and only if
$\lambda'\le\mu'$.

\Proof: First, by applying $\omega^i$ we may assume $i=n$. Let
$N:=\lambda_n=\mu_n$. Then, by applying $\omega^{(n+1)N}$ we may assume
$N=0$. Suppose now that $\lambda,\mu$ is a counterexample. Then, by
applying affine reflections in the first $n-1$ coordinates only and by
using \cite{E42} we may assume that $\mu$ is in the fundamental
alcove, i.e.,
$$
\mu=(\underbrace{x+1,\ldots,x+1}_{a\ \|times|},
\underbrace{x,\ldots,x}_{n-a\ \|times|},0)\quad\hbox{with }0\le a<n.
$$
We necessarily have $d:=|\lambda|=|\mu|=xn+a$. We proceed by induction
on $|d|$, the case $d=0$ being trivial. Assume first that
$x\ge0$. Then there is $j\le n$ with $\lambda_j>0$. After applying the
affine reflection $s_\alpha$ where
$$
\alpha:=\cases{\epsilon_1-\epsilon_j&if $j\le a$ or $a=0$\cr
-\epsilon_1+\epsilon_j+1&otherwise\cr}
$$ to $\lambda$ and $\mu$ we may assume $\lambda_1>0$. That way, we
have
$$
\lambda\le\mu
\Leftrightarrow\omega(\lambda)\le\omega(\mu)\Leftrightarrow
s_n\omega(\lambda)\le s_n\omega(\mu)
{\buildrel{(*)}\over\Leftrightarrow} \omega(\lambda')\le\omega(\mu')
\Leftrightarrow\lambda'\le\mu'
$$
where $(*)$ is the induction hypothesis.

For $x<0$ we proceed similarly. In that case, there is $j\le n$ with
$\lambda_j<0$. Then we use the affine reflection $s_\alpha$ with
$$
\alpha:=\cases{-\epsilon_j+\epsilon_n+1&if $j\le a$\cr
\epsilon_j-\epsilon_n&otherwise\cr}
$$
to obtain $\lambda_n<0$. Finally, we have 
$$
\lambda\le\mu
\Leftrightarrow s_n(\lambda)\le s_n(\mu)\Leftrightarrow
\omega^{-1}s_n(\lambda)\le \omega^{-1} s_n(\mu)
{\buildrel{(*)}\over\Leftrightarrow} \omega^{-1}(\lambda')\le\omega^{-1}(\mu')
\Leftrightarrow\lambda'\le\mu'
$$
\qed

\Corollary inforder. There is a unique order relation on $\Lambda$ whose
restriction to each $\Lambda_n$ is the Bruhat order.

\Proposition kernel. Let $\lambda,\mu\in\Lambda$ with $\lambda\preceq\mu$. Then
$l(\lambda)\ge l(\mu)$.

\Proof: Let $\lambda$, $\mu$ be a counterexample. By
\cite{stableBruhat} we may assume $\lambda_n=0$ and $\mu_n>0$. Then
$(\omega^*)^{-1}(\mu)=(\mu_n-1,\mu_1,\ldots,\mu_{n-1})\in\Lambda$ and
$(\omega^*)^{-1}(\lambda)\preceq(\omega^*)^{-1}(\mu)$. Hence
$(\omega^*)^{-1}(\lambda)\in\Lambda$ by \cite{closed}, i.e.,
$\lambda_n>0$.\qed

\Corollary acc. For every $\lambda\in\Lambda$ there are only finitely many
$\mu\in\Lambda$ with $\lambda\preceq\mu$. In particular, the Bruhat order
on $\Lambda$ satisfies the ascending chain condition.

\Proof: Indeed, $\lambda\preceq\mu$ implies that length and degree of $\mu$
is bounded.\qed

\Corollary dproj. For the Kazhdan\_Lusztig involution holds
$d_{n-1}\pi_n=\pi_nd_n$. In particular, there is an involution $d$ of
$\Mp_\infty$ which is compatible with all $d_n$.

\Proof: \cite{kernel} implies that $d_n$ preserves the kernel of
$\pi_n$. Hence it induces a unique involution $\tilde d$ of
$\Mp_{n-1}$ with $\tilde d\pi_n=\pi_nd_n$. To show $\tilde d= d_{n-1}$
it suffices to show $\tilde d(M_0)=M_0$, $\tilde dH_i=H_i^{-1}\tilde
d$ for $i=1,\ldots,n-2$ and $\tilde d\omega_{n-1}=\omega_{n-1}\tilde
d$. The first statement is clear, the second follows from \cite{E48}:
$$
\tilde dH_i\pi_n=\tilde d\pi_n H_i=\pi_n d_n H_i=\pi_n H_i^{-1} d_n=
H_i^{-1} \pi_n d_n=H_i^{-1}\tilde d\pi_n, 
$$
and the third from \cite{E49}:
$$
\tilde d\omega_{n-1}\pi_n=\tilde d\pi_nH_{n-1}\omega_n=
\pi_n d_n H_{n-1}\omega_n=\pi_nH_{n-1}^{-1}\omega_n d_n=
\omega_{n-1}\pi_n d_n=\omega_{n-1}\tilde d\pi_n.
$$
\qed

\noindent
Let $\Mp_{++}$ be the set of possibly infinite linear combinations
$\sum_{\lambda\in\Lambda}a_\lambda M^\lambda$ with $a_\lambda\in
v\ZZ[v]$.

\Theorem. For every $\lambda\in\Lambda$ there is a unique
$\uM^\lambda\in\Mp_\infty$ with $d(\uM^\lambda)=\uM^\lambda$ and
$\uM^\lambda\in M^\lambda+\cM_{++}$. This element is triangular with
respect to the Bruhat order. Moreover,
$\uM^\lambda=\|lim|_{n\pfeil\infty}\uM^{\lambda_{\le n}}$.

\Proof: For any $n\ge2$ we have
$$57
\pi_n(\uM^{\lambda_{\le n}})=\cases{
\uM^{\lambda_{< n}}&if $\lambda_n=0$\cr0&otherwise\cr}
$$
For $\lambda_0=0$, this follows from \cite{dproj}, otherwise it is
implied by \cite{kernel}. This shows the existence of $\uM^\lambda$
and its triangularity (\cite{inforder}). For uniqueness, suppose there
are two solutions $\uM_1$ and $\uM_2$. Write
$m:=\uM_1-\uM_2=\sum_{\lambda\in\Lambda}a_\lambda M^\lambda$ and let
$\mu$ be maximal with $a_\mu\ne0$ (see \cite{acc}). Then $d(a_\mu)=a_\mu$
and $a_\mu\in v\ZZ[v]$ which is impossible.\qed

An analogous statement holds for Macdonald polynomials:

\Theorem. Let $\lambda\in\Lambda$. Then for any $n\ge2$ we have
$$58
\pi_n(\cE_{\lambda_{\le n}})=\cases{
\cE_{\lambda_{< n}}&if $\lambda_n=0$;\cr0&otherwise.\cr}
$$
In particular, $\cE_\lambda:=\|lim|_{n\pfeil\infty}\cE_{\lambda_{\le
n}}$ exists. Moreover, the recursion formula \cite{E19} is still
valid.

\Proof: Apply $\pi_n$ to both sides of \cite{E19}. If $\lambda_n>0$,
then $m=n$ and $\Phi_m=\overline\Phi_m=\omega$. Thus \cite{E58}
follows from \cite{E50}. Otherwise, we apply \cite{E49}.\qed

For the Cherednik operators we have:

\Proposition.  Let $\xi_i^{(n)}$ be the Cherednik operator \cite{E10} in
$n$ variables. Then the following commutation rules hold:
$$64
\pi_nH_{n-1}\tilde\omega_n=v^{-1}\tilde\omega_{n-1}\pi_n
$$
$$65
\pi_n\xi_i^{(n)}=\xi_i^{(n-1)}\pi_n\quad\hbox{for } i=1,\ldots,n-1
$$
In particular, the limit operator
$\xi_i:=\|lim|_{n\pfeil\infty}\xi_i^{(n)}$ exists and
$$67
\xi_i(\cE_\lambda)=q^{\lambda_i}t^{1-w^\lambda(i)}\cE_\lambda.
$$
Moreover, the $\cE_\lambda$ are, up to a scalar, the only joint
eigenvectors in $\cP_\infty$.

\Proof: By \cite{commdia} we may think of $\pi_n$ as projection
$\cL[z_1,\ldots,z_{n+1}]\pfeil \cL[z_1,\ldots,z_n]$. Then \cite{E18}
shows $\pi_nH_{n-1}=v^{-1}\pi_ns_{n-1}$. A direct calculation using
\cite{E66} shows \cite{E64}. This and the definition \cite{E10} shows
\cite{E65}. Equation \cite{E67} follows readily from
\cite{E13}. Finally assume $\cE$ is another eigenvector. Let
$z^\lambda$ be a monomial occurring in $\cE$ for which $\lambda$ is
maximal with respect to the Bruhat order. The triangularity of $\xi_i$
shows that $\cE$ corresponds to the same eigenvalue as
$\cE_\lambda$. For suitable $a$, the $x^\lambda$\_term of
$\cE':=\cE-a\cE_\lambda$ cancels out. If $\cE'\ne0$ we could replace
$\cE$ by $\cE'$ and obtain a contradiction.\qed

\beginsection almostsymmetric. The almost symmetric submodule

The elements $\uM^\lambda$ cannot form a basis of $\Mp_\infty$ since
that space is far too big. To pin down the span we introduce for any
$\lambda\in\Lambda$ the notation $\lambda_{\le
m}:=(\lambda_1,\ldots,\lambda_m)$ and
$\lambda_{>m}:=(\lambda_{m+1},\lambda_{m+2},\ldots)$. For fixed
$m\ge0$ we define $\cM(m)\subseteq\Mp_\infty$,
$\cP(m)\subseteq\cP_\infty$ as the space of {\it $m$-symmetric
elements}, i.e., elements $\xi$ with
$$
H_i(\xi)=v^{-1}\xi\quad\hbox{for all }i>m.
$$
For $\xi=\sum_{\lambda}a_\lambda M^\lambda\in\Mp_\infty$ this
condition simply means
$$
a_\lambda=v^{\ell(w^\lambda)-\ell(w^\mu)}a_\mu
$$
whenever $\lambda_{\le m}=\mu_{\le m}$ and $\lambda_{>m}$ is a
permutation of $\mu_{>m}$. For $\xi\in\cP_\infty$ it means even simpler that
$\xi$ is symmetric in the variables $z_{m+1},z_{m+2},\ldots$. This
follows from
$$
H_i-v^{-1}=-{v^{-1}z_i-vz_{i+1}\over z_i-z_{i+1}}(1-s_i).
$$

A basis of $\cM(m)$ can be constructed as follows. Let $\Lambda(m)$
be the set of $\lambda\in\Lambda$ such that $\lambda_{>m}$
is a partition. Then, for $\lambda\in\Lambda(m)$ we define
$$
M^{\lambda|m}:=\sum_{\lambda'}v^{\ell(w^{\lambda'})}
M^{\lambda_{\le m}\lambda'},
$$
where $\lambda'$ runs through all permutations of $\lambda_{>m}$ and
where $\lambda_{\le m}\lambda'$ denotes the concatenation
$(\lambda_1,\ldots,\lambda_m,\lambda'_1,\lambda'_2,\ldots)$.

Clearly we have $\cM(0)\subseteq\cM(1)\subseteq\ldots$ and
$\cP(0)\subseteq\cP(1)\subseteq\ldots$. Their unions are denoted by
$\Ma$ and $\Pa$, respectively, and their elements are called ``almost
symmetric''. We still have an isomorphism
$$
\Psi:\Pa\pf\sim \Ma.
$$

Also $\Ma$ possesses a nice basis. For $\lambda\in\Lambda$ we define
its {\it
partition length} $\npd(\lambda)$ as the minimal number $m\ge0$
such that $\lambda_{>m}$ is a partition. For example,
$$
\lambda=(1,2,1,0,2,1,0,0,\ldots)\quad{\rm has}\quad\npd=4.
$$
Moreover, $\npd(\lambda)=0$ if and only if $\lambda$ itself is a
partition. Now we simply define
$$
M^{\lambda|}:=M^{\lambda|\npd(\lambda)}.
$$
For example
$$
M^{(0,2)|}=M^{(0,2)}+vM^{(0,0,2)}+v^2M^{(0,0,0,2)}+\ldots.
$$

\Theorem basis1. The elements $M^{\lambda|}$, $\lambda\in\Lambda$, form an
$\cL$\_basis of $\Ma$. Moreover, the $M^{\lambda|}$ with
$\npd(\lambda)\le m$ span $\cM(m)$.

\Proof: First we show that the $M^{\lambda|}$ span $\Ma$. Clearly, the
$M^{\mu|m}$ with $\mu\in\Lambda(m)$ form a basis of $\cM(m)$. Since
$\Ma$ is the union of the $\cM(m)$ it suffices to show that
$M^{\mu|m}$ is in the span of the $M^{\lambda|}$. If $m=\npd(\mu)$,
then there is nothing to show. Thus assume $m>\npd(\mu)$, i.e.,
$\mu'':=\mu_{\ge m}\in\Lambda(0)$. For each part $a$ of $\mu''$ let
$\mu''_a$ be obtained from $\mu''$ by putting $a$ in front and
omitting one occurrence of $a$. E.g. if $\mu''=(4,3,3,0,\ldots)$, then
$$
\mu''_4=\mu''=(4,3,3,0,\ldots),\
\mu''_3=(3,4,3,0,\ldots),\
\mu''_0=(0,4,3,3,0,\ldots).
$$
Assume $a$ occurs in $\mu''$ for the first time in position
$i_a$. Put $\mu':=\mu_{\le m-1}$ and $\mu_a:=\mu'\mu''_a$.
Then we have the formula
$$
M^{\mu|m-1}=\sum_a v^{i_a-1}M^{\mu_a|m}=M^{\mu|m}+
\sum_{a\ne\mu_m}v^{i_a-1}M^{\mu_a|}
$$
which proves the claim by induction.

As for the linear independence, assume
$\sum_\lambda a_\lambda M^{\lambda|}=0$. Let $\lambda$ be maximal with
respect to the lexicographic order with $a_\lambda\ne0$. Then
$M^\lambda$ occurs only in $M^{\lambda|}$, which yields the
contradiction $a_\lambda=0$.\qed

\noindent The main reason for introducing $\Ma$ is the following

\Theorem. The elements $\uM^\lambda$, $\lambda\in\Lambda$, form an
$\cL$\_basis of $\Ma$. Moreover, the $\uM^\lambda$ with
$\npd(\lambda)\le m$ span $\cM(m)$.

\Proof: \cite{KLsymmetry} implies $\uM^\lambda\in\cM(m)\subseteq\Ma$ for
$m\ge\npd(\lambda)$. Now fix $d\ge0$ and $m\ge0$. Let
$\lambda\in\Lambda(m)$ with $|\lambda|=d$. Then in the expansion
$\uM^\lambda=\sum_\mu a_\mu M^{\mu|}$ only those $\mu$ occur with
$|\mu|=d$, $\mu\in\Lambda(m)$ and $\mu\preceq\lambda$. Moreover,
$a_{\lambda\lambda}=1$. Thus the transition matrix $(a_{\lambda\mu})$
is unitriangular and finite, hence invertible. This implies that every
$M^{\mu|}$ is in the span of the $\uM^\lambda$.\qed

For the Macdonald polynomials we have

\Lemma. The operators $\Phi_m$, and $\overline\Phi_m$
act on $\Pa(n)$ for any $n\ge m$. In particular,
$\cE_\lambda\in\cP_q(n)\subseteq\Pa_q$ for any $n\ge l(\lambda)$.

\Proof: This follows from the fact that $\Phi_m$ and
$\overline\Phi_m$ commute with $H_n$ for any $n>m$.\qed

\noindent Note however that the $\cE_\lambda$ do not span
$\Pa_q$. For example we have

\Lemma. Let $\cP'\subseteq\Pa_q$ be the $\cL_q$\_span of the
$\cE_\lambda$, $\lambda\in\Lambda$. Then $\cP_q(0)\cap\cP'=\cL_q$.

\Proof: Let $\cE=\sum_\lambda c_\lambda\cE_\lambda$ be a finite linear
combination which is not constant. Choose $\lambda\in\Lambda$ with
$c_\lambda$ such that $m:=l(\lambda)$ is maximal. Then $m\ge1$. It is
well known (see, e.g., \cite{Kostka} Thm.~4.2 or \cite{Hmult} below)
that $H_m(\cE_\lambda)=a\cE_{s_m(\lambda)}+b\cE_\lambda$ with
$a\ne0$. Thus, $\cE$ cannot be symmetric.\qed

\beginsection scalarproduct. The scalar product and composition Kostka
functions

Recall the following notation from Macdonald's book \cite{Mac} III.2: for
any integer $m\ge0$ put $\phi_m(t):=\prod_{i=1}^m(1-t^i)$. For a
partition $\lambda\in\Lambda(0)$ and an integer $a\ge0$ let
$m_a(\lambda):=\#\{i\ge1\mid\lambda_i=a\}$ and
$$
b_\lambda(t):=\prod_{a\ge1}\phi_{m_a(\lambda)}(t)
$$

\Theorem. We equip $\QQ(v)$ with the $v$\_adic topology. Then, there is
a unique $\cL$\_linear continuous scalar product
$\Ma\times\Ma\pfeil\QQ(v)$ such that the $M^\lambda$,
$\lambda\in\Lambda$, are orthonormal. It has the property that
$$27
\<M^{\lambda|m},M^{\mu|m}\>={\delta_{\lambda\mu}\over
b_{\lambda_{>m}}(t)}\quad\hbox{for all $m\ge0$ and
$\lambda,\mu\in\Lambda(m)$.}
$$
Moreover, the operators $H_w$, $w\in S_\infty$, are selfadjoint.

\Proof: Uniqueness is clear since the $M^\lambda$ are dense in
$\Ma$. In view of \cite{basis1}, for existence it suffices to show
\cite{E27}. If $\lambda\ne\mu$, then $\lambda_{\le m}\lambda'\ne\mu_{\le
m}\mu'$ where $\lambda'$, $\mu'$ are permutations of $\lambda_{>m}$,
$\mu_{>m}$, respectively. This shows $\<M^{\lambda|m},M^{\mu|m}\>=0$
for $\lambda\ne\mu$ and it remains to compute
$\<M^{\lambda|m},M^{\lambda|m}\>$. For this, we may clearly assume
$m=0$. Then
$$
A_\infty:=\<M^{\lambda|0},M^{\lambda|0}\>=
\sum_{\lambda'\in S_\infty\lambda}v^{2\ell(w^{\lambda'})}.
$$
For $n\ge l(\lambda)$ we let $A_n$ be the subsum with $\lambda'\in
S_n\lambda$. Let $S_\lambda$ be the isotropy group of $\lambda$ in
$S_n$. For a finite Coxeter group $H$ let $p_H(t)$ be the function
$\sum_{w\in H}t^{\ell(w)}$. Then
$$
A_n={p_{S_n}(t)\over p_{S_\lambda}(t)}.
$$
We have
$$
p_{S_n}(t)=(1+t)(1+t+t^2)\ldots(1+t+\ldots+t^{n-1})=
{\prod_{i=1}^n(1-t^i)\over(1-t)^n}.
$$ 
Put $m_a:=m_a(\lambda)$. From $S_\lambda=S_{m_0}\times
S_{m_1}\times\ldots$ with $m_0+m_1+\ldots=n$ we get
$$
p_{S_\lambda}(t)={\prod_{i=1}^{m_0}(1-t^i)\over(1-t)^{m_0}}
{\prod_{i=1}^{m_1}(1-t^i)\over(1-t)^{m_1}}
{\prod_{i=1}^{m_2}(1-t^i)\over(1-t)^{m_2}}\ldots=
{\prod_{i=1}^{m_0}(1-t^i)\over(1-t)^n}b_\lambda(t)
$$
Hence
$$
A_n={\prod_{i={m_0+1}}^n(1-t^i)\over b_\lambda(t)}\Pf{n\pfeil\infty}
{1\over b_\lambda(t)}=A_\infty.
$$

Finally, formula \cite{E25} shows that the matrix of $H_i$ with respect
to the basis $M^\lambda$ is symmetric. This implies that all operators
$H_w$, $w\in S_\infty$, are selfadjoint.\qed

At last, we link Kazhdan\_Lusztig polynomials and Macdonald
polynomials in the following

\Definition: For $\lambda,\mu\in\Lambda$ we define the {\it composition
Kostka function} as
$$
K_{\lambda\mu}(q,t):=\<\uM^\lambda,\Psi(\tilde\cE_\mu)\>.
$$
\medskip

In \cite{Mac0}, Macdonald constructed a two-parameter function
$K_{\lambda\mu}(q,t)$ where $\lambda$ and $\mu$ are partitions and
conjectured that they are polynomials in $q$ and $t$ with
non\_negative integers as coefficients. The fact, that
$K_{\lambda\mu}(q,t)$ is a polynomial was proved almost simultaneously
in \cite{P1}, \cite{P2}, \cite{P3}, \cite{Kostka}, and \cite{P4}. The
remaining positivity conjecture was finally settled affirmatively by
Haiman \cite{Hai}. We are going to show (\cite{Kpart}) that our
$K_{\lambda\mu}$ coincide with Macdonald's in case $\lambda,\mu$ are
partitions. The main ``result'' of this paper is the following

\Conjecture Vermutung2. For all $\lambda,\mu\in\Lambda$ holds
$K_{\lambda\mu}(q,t)\in\NN[v,q]$.

\noindent As for the evidence, we have

\item{$\bullet$} The conjecture is true for $q=0$. In fact,

\Lemma KLpoly. $K_{\lambda\mu}(0,t)$ is a Kazhdan\_Lusztig polynomial.

\Proof: Given $\mu\in\Lambda$, the expansion of the recursion formula
\cite{E21} for $q=0$ gives
$$
\Psi(\tilde\cE_\mu)|_{q=0}=H_{c_{m_d}}\ldots H_{c_{m_2}}H_{c_{m_1}}(M_0)
$$
with uniquely determined numbers $m_d\ge \ldots\ge m_2\ge m_1\ge1$.
It is easy to see that $c_{m_d}c_{m_{d-1}}\ldots c_{m_1}$ is a reduced
decomposition of $m_{-\tau}$. This implies
$$69
\Psi(\tilde\cE_\mu)|_{q=0}=M^\mu
$$
and therefore
$$
\uM^\lambda=\sum_\mu K_{\lambda\mu}(0,t)M^\mu.
$$\qed

\Remark: After the first release of this paper in the arxiv the
specialization statement \cite{E69} has been generalized to arbitrary
root systems by Ion, \cite{Ion}. This led him to speculations about
Macdonald positivity for arbitrary root systems. As explained in the
introduction, such a thing does not even exist for the root system
$A_{n-1}$ where $n$ is fixed. More precisely, any generalization of
Macdonald positivity to arbitrary root systems would require
completely new ideas if it exists at all.

\medskip

\item{$\bullet$} We can almost prove polynomiality. This is our ``real''
main result.

\Theorem. For all $\lambda,\mu\in\Lambda$ holds
$K_{\lambda\mu}(q,t)\in\ZZ[v,v^{-1},q]$.

\Proof: For $m\ge0$ and $\lambda\in\Lambda(m)$ define $\tilde
M^{\lambda|m}:=b_{\lambda_{>m}}(v^2)M^{\lambda|m}$. If $l(\mu)\le m$,
then we can expand
$$
\tilde\cE_\mu=\sum_{\tau\in\Lambda(m)}c_{\mu\tau}\tilde M^{\tau|m}.
$$
It has been shown in \cite{Kostka} Thm.~5.4 that the coefficients
$c_{\mu\tau}$ are in $\ZZ[q,v,v^{-1}]$. If $\lambda$ has
$\npd(\lambda)\le m$, then $\uM^\lambda$ has an expansion in terms of
the $M^{\tau|m}$ with polynomial coefficients (\cite{basis1}). The
claim follows from \cite{E27}.\qed

\Remark: The non\_appearance of negative powers of $v$ is equivalent
to \cite{Vermutung1}. In any case, by tracing through all definitions, it
would be possible to give an explicit upper bound for the pole order of
$K_{\lambda\mu}(q,t)$ at $v=0$ depending only on $\mu$.
\medskip

\item{$\bullet$} Using a computer, we tested the conjecture in
thousands of cases.

\item{$\bullet$} Finally, as mentioned, the conjecture holds for
$\lambda,\mu\in\Lambda(0)$ since Macdonald's Kostka functions are special
cases of ours. We start with a lemma.

\Lemma Hmult. For $i\ge1$, $\mu\in\Lambda$ with $s_i(\mu)\ne\mu$ let
$$
f_\mu:=q^{\mu_i-\mu_{i+1}}t^{w^\mu(i+1)-w^\mu(i)},
\ A_\mu:={v-v^{-1}f_\mu\over1-f_\mu}.
$$
Then
$$63
(H_i-v^{-1})(\cE_\mu)=A_\mu(\cE_{s_i(\mu)}-\cE_\mu).
$$

\Proof: From \cite{Kostka} Theorem~4.2 one deduces the formula
$$62
(H_i-v^{-1})E_{s_i(\mu)}=v^{-1}E_\mu-{v^{-1}-vf_\mu\over1-f_\mu}E_{s_i(\mu)}
$$
under the provision $\mu_i>\mu_{i+1}$. Here $E_\mu$ is the Macdonald
polynomial with the $z^\mu$\_coefficient normalized to $1$. Now
$f_{s_i(\mu)}=f_\mu^{-1}$. Thus replacing $\mu$ by $s_i(\mu)$ in
\cite{E62} results in 
$$61
(H_i-v^{-1})E_{\mu}=v^{-1}E_{s_i(\mu)}-A_\mu E_{\mu}.
$$
Let $c_\lambda$ be the normalization factor \cite{E55}. Then formula
\cite{E63} amounts to $v^{-1}c_\mu/c_{s_i(\mu)}=A_\mu$. This is readily
verified using the fact that $c_\mu$ and $c_{s_i(\mu)}$ differ in only
one factor namely the contribution of the box $(i+1,\mu_i+1)$ and
$(i,\mu_i+1)$, respectively. This proves \cite{E63} in the case
$\mu_i<\mu_{i+1}$. The other case can be easily deduced from that
using the Hecke relation \cite{E1}.\qed

\Corollary Mpart. Let $\lambda,\mu\in\Lambda$, $i\ge1$ with
$\lambda_i\ge\lambda_{i+1}$ and $\mu_i>\mu_{i+1}$. Then
$$
K_{\lambda\,s_i(\mu)}=vK_{\lambda\mu}.
$$

\Proof: Let $\Xi:=H_i-v^{-1}$. Then $\Xi(\uM^\lambda)=0$
(\cite{KLsymmetry}) and
$\Xi(\tilde\cE_\mu)=v^{-1}A_\mu(\tilde\cE_{s_i(\mu)}-v\tilde\cE_\mu)$
(\cite{Hmult}).  Moreover, $\Xi$ is selfadjoint. Hence
$$
0=\<\Xi(\uM^\lambda),\Psi(\tilde\cE_\mu)\>=
v^{-1}A_\mu\<\uM^\lambda,\Psi(\tilde\cE_{s_i(\mu)}-v\tilde\cE_\mu)\>=
v^{-1}A_\mu(K_{\lambda\,s_i(\mu)}-vK_{\lambda\mu}).
$$
\qed

\noindent This result reduces the computation of $K_{\lambda\mu}$ to
the case where $\mu_i\ge\mu_{i+1}$ whenever
$\lambda_i\ge\lambda_{i+1}$. In particular, if $\lambda$ is a
partition, one may assume that $\mu$ is a partition, as well.

Now we introduce the symmetric (i.e., original) Macdonald functions.
Let $\lambda\in\Lambda(0)$ be a partition. The subspace of
$\cL_q[z_1,\ldots,z_n]$ spanned by the $\cE_\mu$, $\mu\in
S_n\lambda_{\le n}$ contains a unique symmetric polynomial
$\cJ_{\lambda_{\le n}}$ whose $z^\lambda$\_coefficient is
$$
\prod_{s\in\lambda}\left(1-q^{a_\lambda(s)}t^{l_\lambda(s)+1}\right).
$$
(note the small difference to \cite{E55}). It follows from \cite{E58}
that the $\cJ_{\lambda_{\le n}}$ are compatible and therefore have a
limit $\cJ_\lambda\in\cP_q(0)$, the {\it symmetric Macdonald function}.

\Lemma. For $\lambda,\mu\in\Lambda(0)$ holds
$\<\uM^\lambda,\Psi(\cJ_\mu)\>=K_{\lambda\mu}$.

\Proof: We work first with a finite number of $n$ variables. Let
$\lambda,\mu\in\Lambda_n$ be partitions. It was shown in the proof of
\cite{Kostka}~Thm.~6.1 that
$$
\cJ_\mu=c_\mu \sum_{w\in W_f}v^{\ell(w)}H_w^{-1}(\cE_{\mu^-})
\quad\hbox{where}
$$
$$
\mu^-:=w_0(\mu),\ \hbox{and}\ 
c_\mu:={(1-t)^n\over
\prod_{i=1}^m(1-q^{\mu_i}t^{n-m+i})\prod_{a=1}^{n-m}(1-t^a)}\
\hbox{with}\ m:=l(\mu).
$$
Thus,
$$
\eqalign{\<\uM^\lambda,\Psi(\cJ_\mu)\>&=c_\mu\sum_{w\in W_f}
v^{\ell(w)}\<H_w^{-1}\uM^\lambda,\Psi(\cE_{\mu^-})\>=\cr
&=d_\mu \<\uM^\lambda,\Psi(\cE_{\mu^-})\>=
d_\mu \<\uM^\lambda,\Psi(\cE_\mu)\>\cr}
$$
where the second equality is \cite{KLsymmetry} and the last is
\cite{Mpart}. The coefficient is
$$
d_\mu=c_\mu\sum_{w\in
W_f}t^{\ell(w)}=c_\mu{\prod_{a=1}^n(1-t^a)\over(1-t)^n}=
\prod_{i=1}^m{1-t^{n-m+i}\over 1-q^{\lambda_i}t^{n-m+i}}.
$$
The assertion follows from $\|lim|\limits_{n\pfeil\infty}d_\mu=1$.\qed

\Theorem Kpart. Assume that both $\lambda$ and $\mu$ are partitions. Then
$K_{\lambda\mu}$ coincides with Macdonald's $q,t$\_Kostka function.

\Proof: On one hand, $\cJ_\mu|_{q=0}$ equals $\tilde M^{\mu|0}$
(\cite{Kostka} Thm.~6.2). On the other hand, it is also the
Hall-Littlewood polynomial $Q_\mu$ (\cite{Mac} VI (8.4)ii). Let
$\<\cdot,\cdot\>_{\rm HL}$ be the scalar product of \cite{Mac}~III.4
on symmetric functions making the Hall-Littlewood functions
orthogonal. Then the comparison of \cite{E27} with the scalar product
of Hall\_Littlewood functions shows
$$
\<\Psi(f),\Psi(g)\>=\<f,g\>_{\rm HL}
$$
for any two symmetric functions $f,g\in\cP(0)$. Since $\lambda\in\Lambda(0)$
we have $\uM^\lambda=\Psi(s_\lambda)$ (\cite{Lusztig}). Thus we have
$$
K_{\lambda\mu}=\<\uM^\lambda,\Psi(\cJ_\mu)\>=\<s_\lambda,\cJ_\mu\>_{\rm HL}.
$$
The last expression is just Macdonald's definition of $K_{\lambda\mu}$.\qed

\beginsection Refinement. A refinement

The recursive formula \cite{E21} can be expanded to give a closed
formula for the polynomials $\tilde\cE_\lambda$. For this we define
the {\it column\_length} of $s\in\lambda$ as
$$
c_\lambda(s):=\#\{k<i\mid j\le\lambda_k+1\}+\#\{k\ge i\mid j\le \lambda_k\}.
$$
If $\lambda$ is a partition, then $c_\lambda(s)$ is the length of the
column containing $s$.  Now we enumerate the boxes of $\lambda$ from
the top to the bottom starting with the {\it right\/}most column and
working to the left. For example $\lambda=(3,0,1,2,0,\ldots)$ gives
$$20
\vcenter{\young(421,,5,63)}
$$
For $i=1,\ldots,|\lambda|$ we put $c_i:=c_\lambda(s_i)$. In the
example above we get the sequence
$$
1,2,3,3,4,4
$$
For
$m\ge1$ we define the operators
$$
X_m^{(0)}:=\Phi_m,\quad X_m^{(1)}:=-\overline\Phi_m.
$$
A {\it marked} diagram $\blambda$ is diagram $\lambda$ together with a
subset $S\subseteq\lambda$ of boxes. Let $\epsilon_i=1$ if $s_i\in S$
and $\epsilon_i=0$ otherwise. Then for a marked diagram
$S\subseteq\lambda$ (with $n:=|\lambda|$) we define the partial
Macdonald polynomial as
$$
\tilde\cE_{\blambda}:=
X_{c_n}^{(\epsilon_n)}\ldots X_{c_1}^{(\epsilon_1)}(1).
$$
For example, the marked diagram of the shape \cite{E20}
$$
\vcenter{\young(\leer\punkt\leer,,\punkt,\leer\punkt)}
$$
gives
$$
\tilde\cE_{\blambda}=
X_4^{(0)}X_4^{(1)}X_3^{(0)}X_3^{(1)}X_2^{(1)}X_1^{(0)}(1)=
-\Phi_4\overline\Phi_4\Phi_3\overline\Phi_3\overline\Phi_2\Phi_1(1).
$$

\Theorem. The Macdonald polynomial $\tilde\cE_\lambda$ can be expressed as
$$60
\tilde\cE_\lambda=\sum_S q^{A_{\blambda}}
t^{L_{\blambda}}\tilde\cE_{\blambda}
$$
where $S$ runs through all markings of $\lambda$ and where
$$
A_{\blambda}:=\sum_{s\in S}(a_\lambda(s)+1),\quad
L_{\blambda}:=\sum_{s\in S}(l_\lambda(s)+1).
$$

\Proof: For a non\_empty diagram $\mu$ we define the following
operation: take the last row, remove its leftmost box and put the
remainder of the row on top, e.g.,
$$
\vcenter{\young(421,,5,63)\qquad\raise15pt\hbox{$\mapsto$}
\qquad\young(3,421,,5)}
$$
The result is $\lambda^*$. Moreover, it is easily verified that the number,
the arm\_length, the leg\_length, and the column\_length of the surviving boxes
don't change. Let $s\in\lambda$ be the bottom left box. Then
$c_\lambda(s)=l(\lambda)=m$, $a_\lambda(s)+1=\lambda_m$, and
$l_\lambda(s)+1=a$ (defined in \cite{E59}). Thus, \cite{E60} is an
expansion of \cite{E21}.\qed

Accordingly, if we define the marked composition Kostka function as
$$
K_{\lambda\bmu}(t):=
t^{L_{\bmu}}\<\uM^\lambda,\Psi(\tilde\cE_{\bmu}),\>,
$$
then we have
$$
K_{\lambda\mu}(q,t)=\sum_S q^{A_{\bmu}}
K_{\lambda\bmu}(t).
$$
The same proof as for $K_{\lambda\mu}$ shows
$K_{\lambda\bmu}\in\ZZ[v,v^{-1}]$. Indeed, the following seems to be true:

\Conjecture. For all $\lambda\in\Lambda$ and all marked diagrams
$\bmu$ holds $K_{\lambda\bmu}\in\NN[v]$.

\noindent Observe that this conjecture is
indeed stronger than \cite{Vermutung2} since there are plenty of
marked diagram with the same exponent $A_{\bmu}$. One of the simplest
examples is $\lambda=(3,1)$ and $\mu=(2,2)$. Here
$$
K_{\lambda\mu}=t+tq+t^2q
$$
and all summands come from different marked diagrams namely
$$
\vcenter{\young(\leer\leer,\leer\leer)\quad 
\young(\leer\punkt,\leer\leer)\quad
\young(\leer\leer,\leer\punkt)}
$$
Another example is
$\lambda=(3,1,1)$, $\mu=(2,2,1)$. Here
$$
K_{\lambda\mu}=t+(t^2+t^3)q+(t+t^2)q+t^3q^2
$$
where the summands come from
$$
\vcenter{\young(\leer\leer,\leer\leer,\leer)\quad
\young(\leer\punkt,\leer\leer,\leer)\quad
\young(\leer\leer,\leer\punkt,\leer)\quad
\young(\leer\punkt,\leer\punkt,\leer)}
$$

\beginrefs

\L|Abk:Ch|Sig:Ch|Au:Cherednik, I.|Tit:Double affine Hecke algebras and
Macdonald's conjectures|Zs:Ann. of Math. (2)|Bd:141|S:191--216|J:1995|xxx:-||

\Pr|Abk:P1|Sig:GR|Au:Garsia, A.; Remmel, J.|Artikel:Plethystic formulas
and positivity for $q,t$-Kostka coefficients|Titel:Mathematical essays
in honor of Gian-Carlo Rota (Cambridge, MA,
1996)|Hgr:-|Reihe:Progr. Math.|Bd:161|Verlag:Birkh\"auser
Verlag|Ort:Boston|S:245--262|J:1998|xxx:-||

\L|Abk:P2|Sig:GT|Au:Garsia, A.; Tesler, G.|Tit:Plethystic formulas for
Macdonald $q,t$-Kostka coefficients|Zs:Adv.
Math.|Bd:123|S:144--222|J:1996|xxx:-||

\L|Abk:Hai|Sig:Ha|Au:Haiman, M.|Tit:Hilbert schemes, polygraphs and the
Macdonald positivity conjecture|Zs:J. Amer. Math.
Soc.|Bd:14|S:941--1006|J:2001|xxx:math.AG/0010246||

\L|Abk:Ion|Sig:Ion|Au:Ion, B.|Tit:A Kato-Lusztig formula for
nonsymmetric Macdonald polynomials|Zs:Preprint|Bd:-|S:28
pages|J:2004|xxx:math.QA/0406060||

\L|Abk:KL|Sig:KL|Au:Kazhdan, D.; Lusztig, G.|Tit:Representations of
Coxeter groups and Hecke algebras|Zs:Invent.
Math.|Bd:53|S:165--184|J:1979|xxx:-||

\L|Abk:P3|Sig:Ki|Au:Kirillov, A.; Noumi, M.|Tit:Affine Hecke algebras and
raising operators for Macdonald polynomials|Zs:Duke Math.
J.|Bd:93|S:1--39|J:1998|xxx:q-alg/9605004||
                                                                        
\L|Abk:Kostka|Sig:\\Kn|Au:Knop, F.|Tit:Integrality of two variable Kostka
functions|Zs:J. Reine Angew. Math.|Bd:482|S:177--189|J:1997|xxx:q-alg/9603027||

\L|Abk:Hecke|Sig:\\Kn|Au:Knop, F.|Tit:On the Kazhdan-Lusztig basis of
a spherical Hecke algebra|Zs:Preprint|Bd:-|S:10
pages|J:2003|xxx:math.RT/0403066||

\L|Abk:Lu1|Sig:\\Lu|Au:Lusztig, G.|Tit:Green polynomials and singularities
 of unipotent classes|Zs:Adv. in Math.|Bd:42|S:169--178|J:1981|xxx:-||

\Pr|Abk:Lu2|Sig:\\Lu|Au:Lusztig, G.|Artikel:Singularities, character
formulas, and a $q$-analog of weight multiplicities|Titel:Analysis
and topology on singular spaces, II, III (Luminy,
1981)|Hgr:-|Reihe:Ast\'erisque|Bd:101-102|Verlag:Soc. Math.
France|Ort:Paris|S:208--229|J:1983|xxx:-||

\L|Abk:Lu3|Sig:\\Lu|Au:Lusztig, G.|Tit:Affine Hecke algebras and their
graded versions|Zs:J. Amer. Math. Soc.|Bd:2|S:599-635|J:1989|xxx:-||

\L|Abk:Mac0|Sig:\\M|Au:Macdonald, I.|Tit:A new class of symmetric
functions|Zs:S\'eminaire lotharingien de
Combinatoire|Bd:B20a|S:41pp.|J:1988|xxx:www.emis.de/journals/SLC/%
opapers/\hfill\break s20macdonald.html||

\B|Abk:Mac|Sig:\\M|Au:Macdonald, I.|Tit:Symmetric functions and Hall
polynomials
(2nd ed.)|Reihe:-|Verlag:Cla\-ren\-don Press|Ort:Oxford|J:1995|xxx:-||

\B|Abk:Mac2|Sig:\\M|Au:Macdonald, I.|Tit:Affine Hecke algebras and
orthogonal polynomials|Reihe:Cambridge Tracts in
Mathematics|Verlag:Cambridge University Press|Ort:Cambridge|J:2003|xxx:-||

\L|Abk:P4|Sig:Sa|Au:Sahi, S.|Tit:Interpolation, integrality, and a
generalization of Macdonald's polynomials|Zs:Internat. Math. Res.
Notices|Bd:10|S:457--471|J:1996|xxx:-||

\L|Abk:Soe|Sig:So|Au:Soergel, W.|Tit:Kazhdan-Lusztig-Polynome und eine
Kombinatorik f\"ur Kipp-Moduln|Zs:Represent.
Theory|Bd:1|S:37--68|J:1997|xxx:-||

\endrefs

\bye